\documentclass[11pt]{amsart}
\usepackage{amsxtra}
\theoremstyle{plain}

\newcommand{\cleqn}{\setcounter{equation}{0}}
\newcommand{\clth}{\setcounter{theorem}{0}}
\newcommand {\sectionnew}[1]{\section{#1}\cleqn\clth}

\newtheorem{theorem}{Theorem}[section]
\newtheorem{lemma}[theorem]{Lemma}
\newtheorem{definition-lemma}[theorem]{Definition-Lemma}
\newtheorem{proposition}[theorem]{Proposition}
\newtheorem{corollary}[theorem]{Corollary}
\newtheorem{definition}[theorem]{Definition}
\newtheorem{example}[theorem]{Example}
\newtheorem{remark}[theorem]{Remark}
\newtheorem{notation}[theorem]{Notation}
\newtheorem{excercise}[theorem]{Excercise}
\newtheorem{assumption}[theorem]{Assumption}
\newcommand \bth[1] { \begin{theorem}\label{t#1} }
\newcommand \ble[1] { \begin{lemma}\label{l#1} }

\newcommand \bpr[1] { \begin{proposition}\label{p#1} }
\newcommand \bco[1] { \begin{corollary}\label{c#1} }
\newcommand \bde[1] { \begin{definition}\label{d#1}\rm }
\newcommand \bex[1] { \begin{example}\label{e#1}\rm }
\newcommand \bre[1] { \begin{remark}\label{r#1}\rm }

\newcommand \bnota[1] { \begin{notation}\label{n#1}\rm }
\newcommand \bexc[1] { \begin{excercise}\label{x#1}\rm }
\newcommand \bas[1] { \begin{assumption}\label{a#1}\rm }
\newcommand {\eth} { \end{theorem} }
\newcommand {\ele} { \end{lemma} }

\newcommand {\epr} { \end{proposition} }
\newcommand {\eco} { \end{corollary} }
\newcommand {\ede} { \end{definition} }
\newcommand {\eex} { \end{example} }
\newcommand {\ere} { \end{remark} }

\newcommand {\enota} { \end{notation} }
\newcommand {\eexc} { \end{excercise} }
\newcommand {\eas} {\end{assumption}}

\newcommand \thref[1]{Theorem \ref{t#1}}
\newcommand \leref[1]{Lemma \ref{l#1}}
\newcommand \prref[1]{Proposition \ref{p#1}}
\newcommand \coref[1]{Corollary \ref{c#1}}

\newcommand \exref[1]{Example \ref{e#1}}
\newcommand \reref[1]{Remark \ref{r#1}}
\newcommand \lb[1]{\label{#1}}

\newcommand \notaref[1]{Notation \ref{n#1}}
\def \d {{\partial}}   

\def \Cset {{\mathbb C}}

\def \C {{\mathbb C}}

\def \tP   {{\wt{P}}}

\def \tG   {{\tilde{G}}}
\def \tB   {{\tilde{B}}}
\def \tT   {{\tilde{H}}}
\def \tU   {{\tilde{U}}}
\def \tg   {{\tilde{g}}}
\def \tt   {{\tilde{t}}}
\def \tu   {{\tilde{u}}}
\def \tBminus   {{\tilde{B}}_-}
\def \tXttw {{\tilde{X}}_{\tt, w}}
\def \tmu {{\tilde{\mu}}}
\def \tXtt {{\tilde{X}}_{\tt}}

\def \tg   {{\tilde{g}}}

\def \calU {{\mathcal{U}}}

\def \calR  {{\mathcal R}}

\def \ra  {\rightarrow}           
\def \lra {\longrightarrow}

\def \wt {\widetilde}

\def \hs {\hspace{.2in}}

\def \Ad { {\mathrm{Ad}} }

\def \g  {\mathfrak{g}}   
\def \h  {\mathfrak{h}}
\def \d  {\mathfrak{d}}

\def \n  {\mathfrak{n}}

\def \b  {\mathfrak{b}}

\def \sl {\mathfrak{sl}}

\DeclareMathOperator \tr { {\mathrm{tr}} }

\def \dw {\dot{w}}

\def \du {\dot{u}}

\def \la {{\langle}}
\def \ra {{\rangle}}
\def \lara {\la \,\, , \, \ra}

\def \ud {\underline{d}}

\def \tG {\tilde{G}}

\def \gdia {\g_{\scriptscriptstyle \Delta}}
\def \gst {\g^*_{{\rm st}}}
\def \Gdia {G_{{\scriptscriptstyle \Delta}}}
\def \piGs {\pi_{{\scriptscriptstyle G}^*}}
\def \piG {\pi_{\scriptscriptstyle G}}
\def \piP {\pi_{\scriptscriptstyle P}}
\def \piR {\pi_{\scriptscriptstyle R}}
\def \piGB {\pi_{\scriptscriptstyle G/B}}
\def \ud {\underline{d}}

\def \tpiP {\tilde{\pi}_{{}_P}}
\def \Bdia {B_{{\scriptscriptstyle \Delta}}}

\def \piX {\pi_{{\scriptscriptstyle X}}}
\def \piY {\pi_{{\scriptscriptstyle Y}}}
\def \piGH {\pi_{{\scriptscriptstyle G/H}}}
\def \Xu {X^u_{t, w_0}}

\begin{document}
\title[Poisson geometry of the Grothendieck resolution]
{Poisson geometry of the Grothendieck resolution of a complex
semisimple group}
\author[Sam Evens, Jiang-Hua Lu]{Sam Evens, Jiang-Hua Lu}
\address{Department of Mathematics, University of Notre Dame, Notre 
Dame, 46556
\newline \indent Department of Mathematics, University of Hong Kong, Pokfulam Road, Hong Kong}

\email{evens.1@nd.edu, jhlu@maths.hku.hk}
\date{}

\begin{abstract}

\noindent
Let $G$ be a complex semi-simple Lie group with a fixed pair of 
opposite Borel subgroups $(B, B_-)$.
We study a Poisson structure $\pi$ on $G$ and a Poisson structure $\Pi$
on the Grothendieck resolution $X$ of
$G$ such that the Grothendieck
map $\mu:(X,\Pi) \to (G,\pi)$ is Poisson. We show that
the orbits of 
symplectic leaves of $\pi$ in $G$ under the conjugation action by 
the Cartan subgroup $H= B \cap B_-$ are  intersections of conjugacy classes
and Bruhat cells $BwB_-$,
while the $H$-orbits of symplectic leaves of $\Pi$ on $X$
give desingularizations of intersections of Steinberg fibers and Bruhat cells in $G$.
We also give birational Poisson isomorphisms from  quotients by $H \times H$ of
products of double 
Bruhat cells in $G$ to intersections of Steinberg fibers and Bruhat cells.
\end{abstract}

\maketitle

{\it Dedicated to Victor Ginzburg for his 50th birthday}

\tableofcontents
\sectionnew{Introduction}\lb{intro}
 
Let $G$ be a complex connected semi-simple Lie group.
It is well-known that
the choice of a pair $(B, B_-)$ of opposite Borel subgroups of $G$ leads to 
a standard Poisson structure $\piG$ on $G$ making $(G, \piG)$ into a 
Poisson Lie group \cite{k-s:quantum}. The Poisson Lie group $(G, \piG)$ is the semi-classical
limit of the quantum group $\Cset_q(G)$, the dual (as Hopf algebras) of the much studied
quantized universal enveloping algebra $U_q \g$, where $\g$ is the Lie algebra of $G$.
In addition to the relationship between the representation theory of
$U_q \g$ and symplectic leaves of $\piG$ established by Hodges,
Levasseur, DeConcini, Kac, Procesi and others (see \cite{HL, DKP}), there are deep
connections between the geometry of $\piG$ and cluster algebras,
total positivity, and canonical bases,
as seen in the fascinating work of Lusztig, Fomin, Zelevinsky
and others
(see, for example, \cite{beren-z, FZ, fz:I, fz:II, lusztig}).
In particular, let $H = B \cap B_-$, a Cartan subgroup of $G$. Then the
$H$-orbits of symplectic leaves of $\piG$ in $G$ are \cite{HL, KZ}
the so-called double Bruhat cells
defined as
\[
G^{u, v} = BuB \cap B_- v B_-, \hs u, v \in W,
\]
where $W$ is the Weyl group of $G$ with respect to $H$.
Double Bruhat cells are the main motivating examples for cluster algebras
\cite{FZ, fz:I, fz:II}. 

Because of the connections between $(G, \piG)$ 
and the quantum groups
and cluster algebras, it is natural to study the
dual Poisson Lie group $(G^*, \pi_{{\scriptscriptstyle G}^*})$ of $(G, \piG)$,
which is the semi-classical limit of $U_q \g$. 
Let $U$ and $U_-$ be the unipotent radicals of $B$ and
$B_-$ respectively. Then $G^*$ is the subgroup of $B \times B_-$ given by
\[
G^* = \{(uh, h^{-1}u_-) \mid u \in U, u_- \in U_-, h \in H\} \subset B \times B_-.
\]
The map $\eta: G^* \to G: \; (b, b_-) \mapsto bb_{-}^{-1}$
is a local diffeomorphism and
$\eta(G^*)=BB_-$ is open in $G$. It turns out that there is a unique 
Poisson structure $\pi$ on $G$
such that $\eta: (G^*, \pi_{{\scriptscriptstyle G}^*}) \to (
G, -\pi)$ is a Poisson map. Thus we may study the Poisson geometry of
$(G^*, \pi_{{\scriptscriptstyle G}^*})$ by studying $(G,-\pi)$.
The first part of
the paper, $\S$\ref{sec-pi0-on-G},
is devoted to the Poisson structure $\pi$ on $G$.
(When $G$ is of adjoint type,
the Poisson structure $\pi$ was considered in our previous papers
\cite{e-l:real} and \cite{e-l:cplx}, where its extension to
the wonderful compactification of $G$ is also discussed).

Recall  that
a {\it Poisson action} of the Poisson Lie group $(G, \piG)$ is an action of $G$ on
a Poisson manifold $(P, \piP)$ such that
the action map $(G, \piG) \times (P, \piP) \to (P, \piP)$ is Poisson. We show in
\prref{pr-G-Gs-poi-on-pi0} and \prref{pr-H-orbits} that

1) the conjugation action of $(G, \piG)$ on $(G, \pi)$ is Poisson. In particular,
$\pi$ is invariant under conjugation by  $H$;

2) the $H$-orbits 
of symplectic leaves of $\pi$ in $G$
are the nonempty intersections of conjugacy classes of $G$
and Bruhat cells $BwB_-$ for $w \in W$.

Since conjugacy classes and Bruhat cells in $G$ do not always intersect
(see \reref{re-not-connected}), we study Steinberg fibers in place
of conjugacy classes.
Intersections of Steinberg fibers and Bruhat
cells in $G$ are studied in $\S$\ref{sec-inter-Bruhat-Steinberg}.
Recall \cite{hum-conjugacy} that a Steinberg fiber in $G$ is the closure of 
a regular conjugacy class in $G$ and is a finite union of
conjugacy classes. It is shown in \prref{pr-inter-Bruhat-Steinberg} that
the intersection of a Steinberg fiber and a Bruhat cell is always nonempty.
For $t \in H$, let $F_t$ be the Steinberg fiber through $t$. Then $F_t = F_{t^\prime}$
if and only if $t^\prime = w_\cdot t$ for some $w \in W$, and 
\[
G = \bigsqcup_{t \in T/W, w \in W} (F_t \cap BwB_-) \hs (\mbox{disjoint union}).
\]
For $t \in H$ and $w \in W$, set
$F_{t, w} = F_t \cap BwB_-$. Since $F_t$ is a finite union of
conjugacy classes, $F_{t, w}$ is a finite union of $H$-orbits of symplectic 
leaves of $\pi$. The subvariety $F_{t, w}$ of $G$ is irreducible 
(\prref{pr-inter-Bruhat-Steinberg}) but may be singular (see \exref{ex-sl3-1}).

Aside from our Poisson geometric motivation, the geometry of Steinberg
fibers is quite important and subtle. In particular, the unipotent
variety is a special case of a Steinberg fiber, and plays a key
role in the Springer correspondence and the study of subregular
singularities \cite{hum-conjugacy, slodowy}.
Further, we believe that intersections  of Steinberg fibers and  
Bruhat cells in $G$ are
important, and that it is worthwhile to study their
algebro-geometric properties as well as the
Poisson structure $\pi$ on them. To do this, we consider in 
$\S$\ref{sec-pi}
the Grothendieck resolution
$X= G \times_B B$ of $G$ with the resolution map 
\[
\mu: \; \; X \longrightarrow G: \; \; \; [g,b] \longmapsto gbg^{-1}, \hs 
g \in G, b \in B,
\]
where $[g b_1, b_1^{-1}bb_1] = [g,b]$ for $g \in G$ and $b, b_1 \in B$.
Note that $\mu$ is $G$-equivariant, where $G$ acts on $G$ by conjugation and on $X$ by
\[
\sigma: \; \; G \times X \lra X: \; \; \; {g_1}_\cdot [g, b] = 
[g_1 g,\,  b], \; \; g_1, g \in G, b \in B.
\] 
We introduce a Poisson structure $\Pi$ on $X$ and show in
\prref{pr-Xtw-smooth-any-G}, \prref{pr-properties-pi}, and \thref{th-H-leaves-Pi}
that $\Pi$ has the following properties:

1) the morphism $\mu:(X,\Pi) \to (G,\pi)$ is Poisson;

2) the action $\sigma$ of $(G, \piG) $ on $(X, \Pi)$ is Poisson.
In particular, $\Pi$ is $H$-invariant;

3) the $H$-orbits of symplectic leaves of $\Pi$ in $X$ are 
the smooth and irreducible subvarieties $X_{t, w}$ of $X$, where
$t \in H, w \in W$, and
\[
X_{t,w} = (G \times_B tU) \cap \mu^{-1}(BwB_-) \subset X.
\]
Assume that $G$ is simply connected.
For $t \in H$, let $X_t = G \times_B tU \subset X$.
It is well-known \cite{slodowy, steinberg-inventiones}
that $\mu: X_t \to F_t$ is a resolution of singularities of $F_t$. We prove in
\coref{co-Xtw-Ftw-resol} that for every $t \in H$ and $w \in W$,
\[
\mu: \; \; X_{t, w} \longrightarrow F_{t, w}
\]
is a resolution of singularities of $F_{t,w}$. Note that 
each $F_{t, w}$ is a {\it finite union}
of $H$-orbits of symplectic leaves of
$\pi$ in $G$, while $X_{t, w}$ is a {\it single} $H$-orbit
of symplectic leaves of $\Pi$ in $X$. Thus, for any $t \in H$ and $w \in W$,
the Poisson morphism 
\[
\mu:  (X_{t, w}, \Pi) \longrightarrow (F_{t, w}, \pi)
\]
may be used to better understand the singular Poisson structure $\pi$ on
 $F_{t,w}$,
and we call it 
 an $H$-equivariant {\it Poisson desingularization} (see
$\S$\ref{poissondesing}). 

In $\S$\ref{sec-log-pi}, the desingularization $\mu:  X_{t, w} \to F_{t, w}$
is used to obtain rational parametrizations of  $F_{t, w}$.
More precisely,
let $w_0$ be the longest element in $W$.  We construct 
an explicit biregular Poisson isomorphism between 
a Zariski open subset of $X_{t, w}$ and the Poisson variety
$(G^{1, w^{-1}w_0}, \, \piG)\times_H (G^{1, w_0},  \, \piG)/H$,
where
$H \times H$ acts on $G^{1, w^{-1}w_0} \times G^{1, w_0}$ from the right by
\[
(g_1, g_2) \cdot (h_1, h_2) = (g_1 h_1, \; h_1^{-1} g_2 h_2), \hs g_1, g_2 \in G, \, h_1,
h_2 \in H.
\]
In  future work, we will use the birational isomorphism to study 
log-canonical coordinates for $(X_{t, w}, \Pi)$ and $(F_{t,w}, \pi)$
and investigate the combinatorial consequences and relations
to work of Fomin, Kogan, and Zelevinsky in \cite{FZ, KZ}.

Our Poisson geometric interpretation of the Grothendieck resolution
$\mu: G \times_B B \to G$ is 
heavily influenced by ideas of Victor Ginzburg, who
emphasized the symplectic nature of the Grothendieck resolution
\[
\mu_0: \; \; G \times_B \b \to \g: \; \; (g,\, x)\longmapsto {\Ad}_g(x), \hs g\in G, 
\,  x \in \g,
\]
in geometry and representation theory \cite{chriss-ginzburg}, 
where  $\g$ and $\b$ are the Lie algebras of $G$ and $B$ respectively.
Note that the vector space $\g$ can be given the linear 
Kostant-Kirillov Poisson structure
$\pi_{0}$,
while $G \times_B \b$ has the Poisson structure $\Pi_0$ whose
symplectic leaves are $G \times_B (y + \n)$ with the 
twisted cotangent bundle symplectic structures,
where $y \in \h$, and $\h$ and $\n$ are the Lie algebras of $H$ and $U$  
respectively. It is easy to see that the linearization of $\pi$ at the identity 
element $e \in G$ is $(\g, \pi_0)$, so $(G, \pi)$ can be
regarded as a deformation of $(\g, \pi_0)$.
In  future work, we hope to study in more detail
the relation between the Poisson morphisms
$\mu: (G \times_B B, \, \Pi) \to (G, \pi)$ and $\mu_0: 
(G \times_B \b, \Pi_0) \to (\g, \pi_0)$.

In $\S$\ref{appendix}, the Appendix,
we collect 
some results on Poisson Lie groups
and coisotropic reduction that are needed in constructing
the Poisson structures $\pi$ on $G$ and $\Pi$ on $G \times_B B$.
We also prove that 
for a Steinberg fiber $F$ and $w \in W$, $F\cap \overline{BwB_-}$ is normal
and Cohen-Macaulay. We compute its dimension and describe its
singular set.

{\bf Acknowledgments:} We would like to thank 
Camille Laurent-Gengoux, Allen Moy, and Milen Yakimov for help and
useful conversations during the preparation of this paper. Both
 authors would like to thank the Hong Kong University of Science
and Technology for its hospitality during the prepartion of this
paper. 
The second author was partially supported by 
HKRGC grants 703304, 703405, and 
the HKU Seed Funding for basic research. 

\subsection{Notation}\lb{intro-notation}
If $G$ is a Lie group with Lie algebra $\g$, and if 
$X \in \wedge^k \g$, $X^L$ and $X^R$ will denote respectively the left
and right invariant $k$-vector fields on $G$ with $X^L(e) = X^R(e) = X$, 
where $e$ is the identity element of $G$.
If $\pi$ is a bi-vector field on a manifold $P$, $\tilde{\pi}$ denotes the bundle map
\begin{equation}\lb{eq-tilde-pi}
\tilde{\pi}: \; \; T^*P \lra TP: \; \; (\tilde{\pi}(\alpha), \; \beta) = \pi(\alpha, \beta), 
\hs \alpha, \beta \in \Omega^1(P).
\end{equation}
By a variety, we mean a complex quasi-projective variety, and a subvariety
is a locally closed subset of a variety. For a morphism $f:X \to Y$
between varieties, we also use   $f:TX \to TY$ for its differential.
When a group $G$ acts on a space $X$,  $g_\cdot x$ will stand for the 
action of $g \in G$ on $x \in X$.

\sectionnew{The Poisson structure $\pi$ on $G$}\lb{sec-pi0-on-G}

Let $G$ be a connected complex semi-simple Lie group. In this section, we  
recall the definition of two Poisson structures $\piG$ and $\pi$ on $G$. The 
Poisson structure $\piG$ on $G$ is {\it multiplicative} \cite{lu-we:poi} and 
$(G, \piG)$ is a Poisson Lie group \cite{k-s:quantum}. 
The Poisson structure $\pi$ has the property that
the conjugation action 
\[
(G, \piG) \times (G, \pi) \lra (G, \pi): \; \; (g, h) \longmapsto ghg^{-1}, \hs g, h \in G
\]
is Poisson. One obtains $\pi$ as a special case of a general construction
in the theory of Poisson Lie groups, which is reviewed
in $\S$\ref{subsec-poi-lie} in the Appendix.

\subsection{The Poisson Lie group $(G, \piG)$ and its dual group 
$(G^*, \pi_{{\scriptscriptstyle G}^*})$}\lb{sec-piG}

Let $\g$ be the Lie algebra of $G$, and let
$\d = \g \oplus \g$ be the direct product Lie algebra. 
Let $\ll \, , \, \gg$ be a fixed non-zero scalar multiple of the 
Killing form of $\g$, and let $\lara$ 
be the symmetric ad-invariant nondegenerate
bilinear form on $\d$ given by
\[
\la x_1 + y_1,  \, \,  x_2 + y_2 \ra = \ll x_1, x_2 \gg - \ll y_1, y_2\gg, 
\hs x_1, x_2, y_1, y_2 \in \g.
\]
Fix a Cartan subalgebra $\h$ of $\g$ and a choice $\Phi^+$ of positive 
roots in the set $\Phi$ of roots for $(\g, \h)$. Let $\g = \h + 
\sum_{\alpha \in \Phi} \g^\alpha$
be the corresponding root decomposition, and let 
\[
\n = \sum_{\alpha \in \Phi^+}\g^\alpha, \hs 
\n_- = \sum_{\alpha \in \Phi^+} \g^{-\alpha}.
\]
The so-called {\it standard Manin triple} associated to $\g$ (see \cite{k-s:quantum})
is the quadruple
$(\d, \gdia, \gst, \lara)$,
where $\g_{{\scriptscriptstyle \Delta}} = \{(x, x): \, x \in \g\}$ is the diagonal of $\d$, and 
\begin{equation}\lb{eq-gstar}
\gst = \h_{-{\scriptscriptstyle \Delta}} + (\n \oplus \n_-) = \{(x + y, -y + x_-): \, 
x \in \n, x_- \in \n_-, y \in \h\}.
\end{equation}
In particular, both $\gdia$ and $\gst$ are maximal isotropic with
respect to $\lara$, and $\lara$  gives rise to a non-degenerate pairing between
$\gdia$ and $\gst$.

Let $G_{{\scriptscriptstyle \Delta}} = \{(g, g): g \in G\} \subset G \times G$, and let
$G^*$ be the connected subgroup of $G \times G$ with Lie algebra $\gst$.
The splitting $\d = \gdia + \gst$ gives rise to multiplicative Poisson structures 
$\piG$ on $G \cong G_{{\scriptscriptstyle \Delta}}$ and $\pi_{{\scriptscriptstyle G}^*}^{}$ on 
$G^*$ making them into a pair of dual Poisson Lie groups \cite{k-s:quantum} (see also the Appendix). 
If
$U$, $U_-$, and $H$ are the connected subgroups of $G$ with Lie algebras $\n, \n_-$, and $\h$
respectively, then 
\begin{equation}\lb{eq-Gstar}
G^* = \{(nh, h^{-1}n_-): \; n \in U, n_- \in U_-, h \in H\} \subset G \times G.
\end{equation}
We will refer to $\piG$ as the {\it standard multiplicative Poisson structure} on $G$.

\bnota{nota-roots}
Throughout the paper, for each $\alpha \in \Phi^+$, we fix root vectors 
$E_{\alpha} \in \g^{\alpha}$ and $E_{-\alpha} \in \g^{-\alpha}$
such that
$\ll E_\alpha, E_{-\alpha} \gg = 1$. 
\enota

The following fact on $\piG$ is
well-known \cite{k-s:quantum}.

\bpr{pr-piG} The Poisson structure $\piG$ is given by 
$\piG = \Lambda_0^R - \Lambda_0^L$, where
$\Lambda_0 = \frac{1}{2} \sum_{\alpha \in \Phi^+} E_\alpha \wedge E_{-\alpha} 
\in \wedge^2 \g$. (See notation in $\S$\ref{intro-notation}).
\epr

\bex{ex-sl2-piG}
Let $G = SL(2, \C) =\left\{\left(\begin{array}{cc} a & b \\ c & d \end{array}\right): 
\; a, b, c, d \in \C, ad-bc = 1\right\},$ and let
$\ll x, y \gg =\frac{1}{2\lambda} \tr(xy)$ for $x, y \in \sl(2, \C)$, where
$\lambda \in \C, \lambda \neq 0$. Let $\h$ be the Cartan subalgebra consisting
of diagonal matrices in $\sl(2, \C)$, and take the standard choices of 
positive and negative roots. Then $\Lambda_0 = \lambda 
\left(\begin{array}{cc} 0 & 1 \\ 0 & 0 \end{array}\right) \wedge 
\left(\begin{array}{cc} 0 & 0 \\ 1 & 0 \end{array}\right)$, and 
the Poisson brackets with respect to
$\piG$ between the functions $a, b, c$ and $d$ are:
\begin{align*}
\{a, b\}_\lambda &= {\lambda}ab, \hs \{a, c\}_\lambda = {\lambda}ac, \; \hs
\{b, d\}_\lambda = {\lambda}bd,\\
\{c, d\}_\lambda &= {\lambda}cd, \hs
\{a, d\}_\lambda = 2{\lambda}bc, \hs \{b, c\}_\lambda = 0.
\end{align*}
We will denote this Poisson structure by $\pi^\lambda_{\scriptscriptstyle SL(2, \C)}$.
\eex

\bnota{Halpha}
For $\alpha \in \Phi^+$, 
let $H_\alpha \in \h$ be such that $\ll H_\alpha, x\gg = \alpha(x)$ for all $x \in \h$, and let
$\ll \alpha, \alpha \gg = \ll H_\alpha, H_\alpha \gg$. Set
\[
h_\alpha = \frac{2}{\ll \alpha, \alpha \gg} H_\alpha, \hs
e_\alpha = \sqrt{\frac{2}{\ll \alpha, \alpha \gg}} E_\alpha, 
\hs e_{-\alpha} = \sqrt{\frac{2}{\ll \alpha, \alpha \gg}} E_{-\alpha}.
\]
Then the linear map 
$\phi_\alpha: \sl(2, \C) \to \g$ given by
\[
\phi_\alpha:   \left(\begin{array}{cc} 1 & 0 \\ 0 & -1 \end{array}\right) \longmapsto h_\alpha, 
\hs 
\left(\begin{array}{cc} 0 & 1 \\ 0 & 0 \end{array}\right) \longmapsto e_\alpha, 
\hs
\left(\begin{array}{cc} 0 & 0 \\ 1 & 0 \end{array}\right) \longmapsto e_{-\alpha}
\]
is a Lie algebra homomorphism. The corresponding Lie group homomorphism from $SL(2, \C)$ to $G$ will
also be denoted by $\phi_\alpha$.
\enota

The following fact is also 
well-known \cite{k-s:quantum}.

\bpr{pr-phi-alpha}
Let $\alpha$ be any simple root, and let $\lambda_\alpha = \frac{\ll \alpha, \alpha \gg}{4}$.
Then the map
\[
\phi_\alpha: \; \; \; \left(SL(2, \C), \; \; \pi^{\lambda_\alpha}_{\scriptscriptstyle SL(2, \C)}
\right)
\lra(G, \, \piG)
\]
is Poisson, where 
$\pi^{\lambda_\alpha}_{\scriptscriptstyle SL(2, \C)}$ is the Poisson structure on $SL(2, \C)$ in
\exref{ex-sl2-piG}.
\epr

Since $\piG$ vanishes at points in $H$, $\piG$ is invariant under left translation by
elements in $H$. By an $H$-orbit of symplectic leaves of $\piG$, we mean a set of the form
$\cup_{h \in H} h\Sigma$, where $\Sigma$ is a symplectic leaf of $\piG$.  For $u, v \in W$, let
$G^{u, v} \subset G$ be the double Bruhat cell given by
\begin{equation}\lb{eq-Guv}
G^{u,v} = BuB \cap B_- v B_-.
\end{equation}
By \cite[Theorem 1.1]{FZ},
$\dim(G^{u,v}) = l(u) + l(v) + \dim(H),$ where $l$ is the length function on $W$.


\ble{le-leaves-piG} \cite{HL, HKKR,  KZ, Reshe}
The $H$-orbits of symplectic leaves of $\piG$ in $G$ are precisely the double Bruhat cells
$G^{u, v}$ for $u, v \in W$.
\ele

\subsection{Definition of the Poisson structure $\pi$ on $G$}\lb{subsec-pi0-on-G}

Let $n = \dim \g$. Let $\{x_j\}_{j = 1}^{n}$ be a basis of $\gdia$ and let
$\{\xi_j\}_{j=1}^{n}$ be the basis of $\gst$ such that $\la x_j, \xi_k\ra = \delta_{jk}$ for 
$1 \leq j, k \leq n$. Let
\begin{equation}\lb{Lambda}
\Lambda = \frac{1}{2} \sum_{j=1}^n (\xi_j \wedge x_j) \in \wedge^2 (\g \oplus \g).
\end{equation}
Let $D = G \times G$. Since $(\g \oplus \g, \gdia, \gst)$ is a Manin triple,
 the bi-vector field
\begin{equation}\lb{eq-piD-plus}
\pi_{\scriptscriptstyle D}^{+} = \Lambda^R + \Lambda^L
\end{equation}
on $D$ is Poisson (see $\S$\ref{subsec-poi-lie} in the Appendix). 
Let $p:D \to D/\Gdia$ be the natural projection. 
Then (see \prref{pr-DtoDG} in the Appendix), 
$p(\pi_{\scriptscriptstyle D}^{+})$
is a Poisson structure on $D/\Gdia$.
Since $p (\Lambda^L) = 0$, 
we also have
$p(\pi_{\scriptscriptstyle D}^{+}) =p(\Lambda^R).$

\bnota{nota-pi0}
Identify $D/\Gdia \cong G$ via the map
$(g_1, g_2) \Gdia \mapsto g_1 g_{2}^{-1}$, and let
\begin{equation}\lb{eq-eta}
\eta: \; \; D=G \times G \longrightarrow G: \; \; (g_1, g_2) \longmapsto g_1 g_{2}^{-1}.
\end{equation}
Througout this paper, $\pi$ will denote the Poisson structure 
$\eta(\pi_{\scriptscriptstyle D}^{+})$  on $G$.
\enota

\bre{Gs}
By \leref{le-Gs-pi0} in the Appendix, the local diffeomorphism
\[
\eta|_{G^*}: \; \; \; (G^*, \, -\piGs) \longrightarrow (G, \pi): \; \; (b, b_-) \longmapsto
bb_-^{-1}
\]
is Poisson.
\ere

Let $\{y_i\}_{i = 1}^{r}$ be a basis of $\h$ such that
$2\ll y_i, y_j\gg = \delta_{ij}$ for $1 \leq i, j \leq r = \dim \h$. As bases of
$\gdia$ and $\gst$, let
\begin{align}\lb{eq-x-i}
\{x_i\}& = \{(y_1, y_1),  \;(y_2, y_2), \ldots, (y_r, y_r),  \;
(E_\alpha, E_\alpha), \;(E_{-\alpha}, E_{-\alpha}): 
\alpha \in \Phi^+\}\\
\lb{eq-xi-i}
\{\xi_i\} & = \{(y_1, -y_1),  \;(y_2, -y_2), \ldots, (y_r, -y_r),\;  
(0, -E_{-\alpha)},\; (E_{\alpha}, 0): 
\alpha \in \Phi^+\}.
\end{align}
The element $\Lambda \in \wedge^2 (\g \oplus \g)$ in \eqref{Lambda}
is then given by
\begin{align}\lb{eq-Lambda}
\Lambda &= \frac{1}{2} \sum_{i=1}^{r} (y_i, -y_i) \wedge (y_i, y_i) \\ \nonumber
\lb{eq-Lambda}& \; \; + 
\frac{1}{2} \sum_{\alpha \in \Phi^+} \left( (E_\alpha, 0) \wedge (E_{-\alpha}, E_{-\alpha}) 
+ (0, -E_{-\alpha}) \wedge (E_\alpha, E_\alpha)\right).
\end{align}
Since $\eta((x,0)^R)=x^R$ and $\eta((0,x)^R)=-x^L$ for $x \in \g$, 
\begin{equation}\lb{eq-pi0-explicit}
\pi = \sum_{i=1}^{r} y_i^L \wedge y_i^R + \frac{1}{2} \sum_{\alpha \in \Phi^+}
(E_\alpha^R \wedge E_{-\alpha}^R + E_\alpha^L \wedge E_{-\alpha}^L) 
+ \sum_{\alpha \in\Phi^+} E_{-\alpha}^{L} \wedge E_\alpha^R.
\end{equation}

\bex{ex-sl2-pi0}
Let $G = SL(2, \C)$. Using $\ll x, y \gg = \tr(xy)$ for $x, y \in \sl(2, \C)$, we can
compute the Poisson structure $\pi$ on
$SL(2, \C)$ directly from \eqref{eq-pi0-explicit} to get 
\begin{align*}
\{a, b\}_\pi &= bd, \hs \{a, c\}_\pi = -cd, \hs \{a, d\}_\pi = 0,\\
\{b, c\}_\pi &= ad-d^2, \hs \{b, d\}_\pi = bd, \hs \{c, d\}_\pi = -cd.
\end{align*}
It is easy to see that $a+d$ is a Casimir function. 
\eex

The following \prref{pr-G-Gs-poi-on-pi0} is a direct consequence of
\prref{pr-DtoDG} in the Appendix and the fact that 
$\piG$  vanishes at every point in $H$.

\bpr{pr-G-Gs-poi-on-pi0}
1) The following group actions are Poisson :
\begin{align}
\lb{G-on-G}
&(G, \piG) \times (G, \pi) \lra (G, \pi): \; \; \; 
(g_1, g) \longmapsto g_1 g g_{1}^{-1},\\
\lb{Gs-on-G}
& (G^*, -\piGs) \times (G, \pi) \lra (G, \pi): \; \; \; 
((b, b_-), g) \longmapsto b g b_{-}^{-1};
\end{align}

2) The Poisson structure $\pi$ on $G$ is invariant under conjugation by $H$. 
\epr

Recall the definition of a coisotropic submanifold
of a Poisson manifold from $\S$\ref{subsec-coiso-reduction} in the Appendix.
The following property of $\pi$ will be used in $\S$\ref{subsec-de-pi}.

\ble{le-tN-coisotropic}
For every $t \in H$, $tU$ is a coisotropic submanifold of $(G, \pi)$.
\ele

\begin{proof}
By \reref{re-coisotropic} in the Appendix, it suffices to show that $\pi(b) 
\in T_{b}(tU) \wedge T_bG$ for every $b = tn \in tU$ with $n \in U$. 
First note that for $1 \leq i \leq r$, 
\[
y_i^L(b) \wedge y_i^R(b) = y_i^L(b) \wedge (y_i^R - y_i^L)(b) 
\in T_b(tU) \wedge T_b G.
\]
Now since $E_\alpha^R(b) \in T_b(tU)$ and
$E_{\alpha}^L(b) \in T_b(tU)$ for any $\alpha \in \Phi^+$, it
follows from 
\eqref{eq-pi0-explicit}
that  $\pi(b) \in   T_b(tU) \wedge T_b G$.
\end{proof}



\subsection{Symplectic leaves of $\pi$ in $G$}\lb{sec-leaves-pi0}

The following \prref{pr-leaves-pi0-1} is a direct consequence 
of \prref{pr-DtoDG} of the 
Appendix. See also \cite{anton-malkin, lu-yakimov-3}.

\bpr{pr-leaves-pi0-1}
The symplectic leaves of $\pi$ in $G$ are precisely the 
connected components of intersections
of conjugacy classes in $G$ and the $G^*$-orbits in $G$, 
where $G^*$ acts on $G$ by \eqref{Gs-on-G}.
\epr

The $G^*$-orbits in $G$
can be determined by the Bruhat decomposition of $G$. Indeed,
let $W$ be the Weyl group of $(G, H)$. For each $w \in W$, 
fix a representative $\dw$ of $w$ in the normalizer
$N_G(H)$ of $H$ in $G$, and let
\begin{equation}\lb{eq-Hw}
U^w = U \cap \dw U\dw^{-1} \subset U \hs \mbox{and} \hs
H_w = \{h w(h): h \in H\} \subset H.
\end{equation}

\ble{le-Gs-orbits-on-G} 1) Every $G^*$-orbit in $G$ for the action 
in \eqref{Gs-on-G} is of the form
$G^*{}_\cdot (h\dw)$ for a unique $w \in W$ and for some $h \in H$. Moreover, 
for any $w \in W$ and $h_1, h_2 \in H$,
$G^*{}_\cdot (h_1 \dw) = 
G^*{}_\cdot (h_2 \dw)$ if and only if $h_1h^{-1}_2 \in H_w$;

2) For any $w \in W$ and $h \in H$, the map
\begin{equation}\lb{eq-phiw-h}
\phi_{w, h}: \; \; U^w \times H_w \times U_- \lra G^*{}_\cdot (h\dw):\; \; 
(n, h_1, n_-) \longmapsto n h h_1 \dw n_-^{-1}
\end{equation}
is a biregular isomorphism.
\ele

\begin{proof}
By the  Bruhat decomposition,
$G = \sqcup_{w \in W} BwB_-$ as a disjoint union.
Clearly each $BwB_-$ is a union of $G^*$-orbits. Moreover, for each $w \in W$, 
\begin{equation}\lb{eq-phiw}
\phi_w:  \; U^w \times H \times U_- \lra BwB_-:   \; (n, h, n_-) \longmapsto
n h \dw n_-
\end{equation}
is a biregular isomorphism.
Since $U\times U_- \subset
G^*$, every $G^*$-orbit in $BwB_-$ is of the form 
$G^*{}_\cdot (h\dw)$ for some $h \in H$. Now let $h_1, h_2 \in H$. Then 
$G^*{}_\cdot (h_1\dw)=G^*{}_\cdot (h_2\dw)$ if and only if there exist $n \in U, n_- \in 
U_-$ and $h \in H$ such that $h_1\dw = nhh_2 \dw h n_-$ which is equivalent to 
$h_1 h_2^{-1} \in H_w$. 
2) follows from the fact that $\phi_w$ is a biregular isomorphism.
\end{proof}

\bre{re-not-connected} A conjugacy class in $G$ may not intersect 
every $G^*$-orbit in $G$.
For example, the conjugacy class of the identity
element of $G$ intersects $G^*{}_\cdot (h\dw)$ if and only if $w = 1$. 
Moreover, the intersection of a conjugacy class and a $G^*$-orbit
in $G$ may not be connected. As an example, consider $G = SL(3, \C)$, and 
let $C$ be the conjugacy class of subregular unipotent elements, i.e., 
\[
C = \{ g \in SL(3, \C): \; \; (g-I_3)^2 = 0, \; g \neq I_3\},
\]
where $I_3$ is the $3 \times 3$ identity matrix. Let $B$ and $B_-$
be the subgroups of all upper and lower
triangular matrices in $SL(3, \C)$ respectively, and let $H = B \cap B_-$.
Let $w_0$ be the longest element in the Weyl group $W \cong S_3$ and let
\[
\dw_0 = \left(\begin{array}{ccc} 0 & 0 & 1 \\ 0 & 1 & 0 \\ -1 & 0 & 0\end{array}
\right).
\]
Fix 
$h = {\rm diag}(h_1, h_2, h_2) \in H$. It is easy to see that 
\[
G^*{}_\cdot (h\dw_0) = \left\{ g = \left(\begin{array}{ccc} a & b & h_1 x\\
 c & h_2 x^{-2} & 0 \\ -h_3 x & 0 & 0 \end{array}\right): 
 \; \; a, b, c, x \in \C, x \neq 0\right\}.
\]
One then checks
 that an element $g \in G^*{}_\cdot (h\dw_0)$ lies in $C$ if and only if
$a = 2, b = c = 0$ and $x^2 = h_2$. Thus the intersection 
$C \cap (G^*{}_\cdot (h\dw_0))$ consists
of exactly two points. 
\ere

\subsection{Intersections of conjugacy classes and Bruhat cells}\lb{subsec-C-BwB}

In contrast to the situation for $G^*$-orbits in $G$,  the intersection 
of a conjugacy class and a Bruhat cell
$BwB_-$ is always connected, as stated in the following \prref{pr-smoothirredint}.

\bpr{pr-smoothirredint}
Let $C$ be a conjugacy class in $G$ and let $w\in W$. Assume that
$C \cap (BwB_-)$ is nonempty.
Then $C \cap (BwB_-)$ is smooth and irreducible, and
$\dim(C \cap (BwB_-)) = \dim(C) - l(w)$.
\epr

\begin{proof}
This follows from \cite[Corollary 1.5]{Ri}
because $G_{{\scriptscriptstyle \Delta}} \cap (B \times B_-)$ is connected. See also
\cite[Proposition 4.10]{e-l:cplx}. 
\end{proof}

Recall \cite{hum-conjugacy} that a conjugacy class $C$ in $G$ is said to be
regular 
if $\dim C = \dim G - \dim H$.

\bpr{pr-regular-C-intersects-Bruhat}
If $C$ is a regular conjugacy class in $G$, then  
$C \cap (BwB_-) \neq \emptyset$ for every $w \in W$.
\epr

\begin{proof}
By \cite[Proposition 5.1 ]{ellers-gordeev}, $C \cap (BwB)\neq \emptyset$.
It follows that $C \cap (Bw)\neq \emptyset$, which implies the result.
\end{proof}

\subsection{$H$-orbits of symplectic leaves of $\pi$}\lb{subsec-H-orbits-pi0}

Since the Poisson structure $\pi$ is invariant under conjugation by elements in $H$, 
if $\Sigma$ is a symplectic leaf of $\pi$, so is $h\Sigma h^{-1}$ for every $h \in H$. 
Let $\Sigma$ be a symplectic leaf of $\pi$. The set 
\[
H {}_\cdot \Sigma := \{h\Sigma h^{-1}: \,h \in H\}
\]
will be called an {\it $H$-orbit of symplectic leaves} of $\pi$ in $G$.
By \prref{pr-leaves-pi0-1} and \leref{le-Gs-orbits-on-G},
 every $H$-orbit of symplectic leaves of $\pi$ in $G$
is contained in $C \cap (BwB_-)$ for some conjugacy class $C$ in $G$ and
some $w \in W$. When $G$ has trivial center, the following \prref{pr-H-orbits} 
is a special case of 
\cite[Corollary 4.7 and Theorem 4.14]{e-l:cplx}. 

\bpr{pr-H-orbits} 
Let $C$ be a conjugacy class in $G$ and let $w\in W$ be such that 
$C \cap (BwB_-) \neq \emptyset$. Then

1) every symplectic leaf of $\pi$ in $C \cap (BwB_-)$ has dimension
equal to 
\[
\dim (C \cap (BwB_-)) - \dim(H/H_w) = \dim C -l(w) -\dim(H/H_w);
\]

2) $C \cap (BwB_-)$ is a single $H$-orbit of symplectic leaves of
$\pi$ in $G$.
\epr

\begin{proof} 1) Let $\Sigma$ be a symplectic leaf of $\pi$ in $C \cap (BwB_-)$.
By \prref{pr-leaves-pi0-1}, $\Sigma$ is a connected component of $C \cap (G^*{}_\cdot (h\dw))$ 
for some $h \in H$. Since $C$ and $G^*{}_\cdot (h\dw)$ intersect transversally, 
using \leref{le-Gs-orbits-on-G} and 
\prref{pr-smoothirredint}, one gets
\begin{align*}
\dim (\Sigma) &= \dim(C \cap (G^*{}_\cdot (h\dw)) = 
\dim(C) + \dim(G^*{}_\cdot (h\dw)) - \dim(G)\\
&=\dim(C) - l(w) - \dim(H/H_w)= \dim (C \cap (BwB_-)) - \dim(H/H_w).
\end{align*}

2)
It is clear that if $\Sigma$ and $\Sigma^\prime$ are two symplectic leaves of
$\pi$ in $C \cap (BwB_-)$, then  either
$H{}_\cdot \Sigma = H {}_\cdot \Sigma^\prime$ or 
$(H{}_\cdot \Sigma) \cap  (H {}_\cdot \Sigma^\prime)
=\emptyset$. Since $C \cap (BwB_-)$ is connected by \prref{pr-smoothirredint}, to prove 2),
it suffices to show that $H {}_\cdot \Sigma$ is open in $C \cap (BwB_-)$
for every symplectic
leaf $\Sigma$ in $C \cap (BwB_-)$. 
To this end, we show that 
the map $H \times \Sigma \to C \cap (BwB_-)$ for the conjugation action by $H$ 
is a submersion. 

Consider the action map
$\alpha:  H \times BwB_- \to BwB_- $ of $H$ on $BwB_-$ by conjugation. 
For $g \in BwB_-$, let 
$\alpha_*: \h \times T_g(BwB_-) 
\to T_g (BwB-)$ be the differential of $\alpha$ at $(e, g) \in H \times BwB_-$,
and let $H {}_\cdot g$ be the $H$-orbit in $BwB_-$ through $g$. Then for any submanifold
$S$ of $BwB_-$ and $g \in S$,
$\alpha_*(\h \times T_g S) = T_g(H {}_\cdot g) + T_g S
\subset T_g(BwB_-),$ so
\begin{equation}\lb{dim}
\dim  \alpha_*(\h \times T_g S) = \dim (T_g(H {}_\cdot g)) + \dim (T_g S) - \dim
(T_g(H {}_\cdot g) \cap  T_g S). 
\end{equation}
First consider the case when $S = G^*{}_\cdot (h \dw)$ for some $h \in H$. 
By using the
isomorphism $\phi_w: U^w \times H \times U_- \to BwB_-$ in \eqref{eq-phiw}, the map
$\alpha: H \times BwB_- \to BwB_-$ becomes the map 
$\beta:  H \times (U^w \times H \times U_-) \to U^w \times H \times U_-$ given by
\begin{equation}\lb{beta}
\beta(h_1, (n, h_2, n_-)) = (h_1 n h_1^{-1}, \, h_2 h_1 w(h_1^{-1}), \, h_1 n_- h_1^{-1}),
\end{equation}
where $h_1, h_2 \in H, n \in U^w$ and $n_- \in U_-.$ 
Let $\h_w = \{x + w(x): x \in \g\}$ be the Lie algebra of $H_w$ and let
$\h_w^- = \{x - w(x): x \in \h\}$. By using the isomorphism
$\phi_{w, h}: U^w \times H_w \times U_- \to G^*{}_\cdot (h \dw)$
and the fact that 
$\h = \h_w + \h_w^-$, one sees from \eqref{beta}  that
$\alpha_*(\h \times T_g (G^*{}_\cdot(h\dw))) = T_g(BwB_-)$ for every 
$g \in  G^*{}_\cdot(h\dw)$. It follows from 
\eqref{dim} that for every $g \in G^*{}_\cdot(h\dw)$,
\begin{align*}
\dim(T_g(H {}_\cdot g) \cap  T_g (G^*{}_\cdot(h\dw))) &=\dim (T_g(H {}_\cdot g)) + \dim (T_g 
(G^*{}_\cdot(h\dw))) - \dim BwB_- \\
&= \dim (T_g(H {}_\cdot g)) - \dim(H/H_w).
\end{align*}
Now let $S = \Sigma$, a connected component of $C \cap (G^*{}_\cdot(h\dw))$. Then by
\eqref{dim}, 
\begin{align*}
\dim \alpha_*(\h \times T_g \Sigma) &\geq \dim (T_g(H {}_\cdot g)) + \dim (T_g \Sigma) - 
\dim(T_g(H {}_\cdot g) \cap  T_g (G^*{}_\cdot(h\dw))) \\
& = \dim (T_g \Sigma) +\dim(H/H_w)= \dim (C \cap (BwB_-))
\end{align*}
for every $g \in \Sigma$, where in the last step we used 1).
Since $\alpha(H \times \Sigma) \subset C \cap (BwB_-)$, we obtain
$\dim \alpha_*(\h \times T_g \Sigma) = \dim (C \cap (BwB_-))$. Thus, the
map
$\alpha|_{H \times \Sigma}: H \times \Sigma \to C \cap (BwB_-)$ has surjective differential at $(e, g)$ 
for every $g \in \Sigma$. It follows 
that the differential of $\alpha|_{H \times \Sigma}$ is surjective
everywhere in $H \times \Sigma$, so
 $H {}_\cdot \Sigma$ is open in $C \cap (BwB_-)$.
This finishes the proof of 2). 

\end{proof}

\bco{R}
The intersection of a regular conjugacy class and a $G^*$-orbit in $G$ is nonempty.
\eco

\begin{proof} Let $C$ be a regular conjugacy class and let 
$G^*{}_\cdot (h\dw)$ be a $G^*$-orbit
in $G$, where $h \in H$ and $w \in W$.  By 
\prref{pr-regular-C-intersects-Bruhat},
$C \cap (BwB_-) \neq \emptyset$. By \leref{le-Gs-orbits-on-G}
and Equation \eqref{beta} in the proof of \prref{pr-H-orbits}, 
$H {}_\cdot (G^*{}_\cdot (h\dw)) = BwB_-$. It follows that 
$C \cap (G^*{}_\cdot (h\dw)) \neq \emptyset$.
\end{proof}

\sectionnew{Intersections of Steinberg fibers and Bruhat 
cells}\lb{sec-inter-Bruhat-Steinberg}

We now apply our results on conjugacy classes to the study of Steinberg fibers.
We will see in this section that
the intersection $F \cap (BwB_-)$ of a Steinberg fiber $F$ and a Bruhat cell 
$BwB_-$ is
always a nonempty irreducible 
Poisson subvariety of $(G, \pi)$. The remainder of the
paper will be devoted to the study of the algebro-geometric and 
Poisson geometric properties of the Poisson varieties $(F \cap (BwB_-), \, \pi)$. 
We do so by
constructing in $\S$\ref{sec-pi} an $H$-equivariant Poisson desingularization of
$(F \cap (BwB_-), \, \pi)$ and in $\S$\ref{sec-log-pi} a birational Poisson isomorphism 
between $(F \cap (BwB_-), \, \pi)$ and 
the $(H \times H)$-quotient of a product of double Bruhat cells in $G$.

\subsection{Steinberg fibers}\lb{subsec-steinberg-fibers}
Recall that a regular class function on $G$ is a regular function on $G$ that is
invariant under conjugation. Two elements $g_1, g_2 \in G$ are said to be in the
same Steinberg fiber \cite[Section 6]{steinberg-IHES} if $f(g_1) = f(g_2)$ for
every regular class function $f$ on $G$.

\bpr{pr-Fz} \cite[6.11 and 6.15]{steinberg-IHES}
Let $F$ be a Steinberg fiber in $G$. Then

1) $F$ is a closed irreducible subvariety of $G$ with codimension $r=\dim H$;

2) $F$ is a finite union of conjugacy classes;

3) $F$ contains a single regular conjugacy class which is dense and open in $F$.
\epr

We now state a lemma that will be used several times in the paper.

\ble{le-irred}
Let Y be an algebraic variety. Assume that $k\geq 0$ is an integer such that
1) each irreducible component of $Y$ has dimension at least $k$, and 2)
$Y = Y_1 \cup Y_2 ... \cup Y_n$ is a disjoint union of subvarieties,
where  $Y_1$ is irreducible and 
$\dim Y_i < k$ for $i \ge 2$.
Then $Y=\overline{{Y}_1}$ is irreducible.
\ele

\bpr{pr-inter-Bruhat-Steinberg}  For any Steinberg fiber $F$ in $G$ and
$w \in W$, 

1) $F \cap (BwB_-)$ is a nonempty irreducible subvariety of $G$ with dimension equal to
$\dim G - \dim H - l(w)$;

2) $F \cap (BwB_-)$ is a finite union of $H$-orbits of symplectic leaves of $\pi$ in $G$.
\epr

\begin{proof} By \prref{pr-regular-C-intersects-Bruhat}
and 3) of \prref{pr-Fz}, $F \cap (BwB_-)$ is nonempty.
Since $\dim BwB_- = \dim G - l(w)$, each irreducible component of
$F \cap BwB_-$ has dimension no less than $\dim G - \dim H - l(w)$
using \prref{pr-Fz}. Decompose $F \cap BwB_- = \cup_{i=1}^n (C_i \cap BwB_-)$,
where $C_1, \dots, C_n$ are the conjugacy classes in $F$. Part 1) follows
by \prref{pr-smoothirredint} and \leref{le-irred}.
Part 2) follows from \prref{pr-H-orbits}, since
 $F$ is a finite union of conjugacy classes.
\end{proof}
 
\subsection{The varieties $F_{t, w}$ and $X_{t, w}$}\lb{subset=Ftw-Xtw}

For $t \in H$, let $F_t$ be the Steinberg fiber containing $t$. 
By the Jordan decomposition
of elements in $G$, every Steinberg fiber is of the form 
$F_t$ for some $t \in H$, and
$F_{w(t)} = F_t$ for $t \in H$ and $w \in W$.

Let $B$ act on $G \times B$ from the right by
\begin{equation}\lb{eq-B-on-GB}
(G \times B) \times B \lra G \times B: \; \; \; ((g, b), \,b_1) 
\longmapsto (gb_1, \, b_1^{-1} b b_1). 
\end{equation}
Denote the $B$-orbit through the point $(g, b) \in G \times B$ by $[g, b]$. 
The map
\begin{equation}\lb{eq-mu}
\mu: \; \; \; G \times_B B \lra G: \; \; \; [g, b] \longmapsto gbg^{-1}.
\end{equation}
is called the {\it Grothendieck (simultaneous) resolution} of $G$.
Since $G$ is the union of all Borel subgroups of $G$ \cite[Page 69]{steinberg-verlag},
$\mu$ is surjective. Note that
$\mu$ is $G$-equivariant, where $G$ acts on $G$ by conjugation and on
$G \times_B B$ by
\begin{equation}\lb{eq-sigma}
\sigma: \; \; G \times (G \times_B B) \longrightarrow G \times_B B: \; \; \; g_1 {}_\cdot [g, b]) =
[g_1g, b], \hs g, g_1 \in G, b \in B.
\end{equation}

For $t \in H$, since $tU$ is invariant under conjugation by $B$, $G \times_B tU$ 
is a smooth subvariety 
of $G \times_B B$. Set
\[
X_t:= G \times_B tU. 
\]
Then $\mu|_{X_t}: X_t \to F_t$ is onto.
Clearly, $G \times_B B = \bigcup_{t \in H} X_t$ is a disjoint union. 

\bnota{nota-Ftw-Xtw}
For $t \in H$ and $w \in W$, let
\[
F_{t, w} = F_t \cap BwB_- \subset G \hs \mbox{and} \hs 
X_{t, w} = X_t \cap \mu^{-1}(BwB_-) \subset G \times_B B.
\]
\enota
Note that for any $t \in H$ and $w \in W$, $X_{t, w} \neq \emptyset$  because
$F_{t, w} \neq \emptyset$.

\subsection{The singularities of $F_{t, w}$}\lb{singular-Ftw}

In this subsection, we assume that $G$ is simply connected.

A Steinberg fiber $F_t$, where $t \in H$,
may be singular, and it is shown in \cite[4.24]{hum-conjugacy}
that the set of smooth points in $F_t$ is  the unique regular 
conjugacy class contained in $F_t$.
The following \prref{pr-smooth-points-FzBwB} concerning the singularities of 
$F_{t,w} = F_t \cap (BwB_-)$ follows from 
\thref{th-completeintersection}, \prref{pr-smoothlocus}
and \thref{th-normality} in the Appendix.

\bpr{pr-smooth-points-FzBwB}  
For any $t \in H$ and $w \in W$,
the set of smooth points of $F_t \cap BwB_-$ is  
$R_t \cap (BwB_-)$, where $R_t$ is the unique regular conjugacy class
 in $F_t$. Moreover,
 $F_t\cap (BwB_-)$ is normal and is a complete
intersection in $BwB_-$.
\epr

\bex{ex-sl3-1}
Consider again $G = SL(3, \C)$ with $B$ and $B_-$ being the subgroups of upper
and lower triangular matrices respectively. 
Let $w_0$ be the longest element in the Weyl group $W$. 
Then
\[
Bw_0B_- = Bw_0 = \left\{\left(\begin{array}{ccc} a & b & y \\ c & x & 0 \\
-\frac{1}{xy} & 0 & 0 \end{array}\right): \; a, b, c, x, y \in \C, x \neq 0, y \neq 0\right\}.
\]
Let $\calU$ be the unipotent subvariety of $G$, so $\calU$ is the Steinberg fiber
containing the identity. Let $C_r$ and $C_{{\rm sr}}$ be the regular and the sub-regular
conjugacy class in $\calU$ respectively. Then ${\calU} \cap (Bw_0)=
(C_r \cap (Bw_0)) \cup (C_{{\rm sr}} \cap (Bw_0))$ and ${\calU} \cap (Bw_0)$
can be identified with the
subset of ${\mathbb C}^5$ with coordinates $(a, b, c, x, y)$ given by
\[
\begin{cases} a + x = 3,\\ ax -bc + \frac{1}{x} = 3, \\ x \neq 0, y \neq 0.
\end{cases}
\]
Using the Jacobian criterion, one sees that
 ${\calU} \cap (Bw_0)$ is singular exactly at the subregular elements
\[
C_{{\rm sr}} \cap (Bw_0) = \left\{\left(\begin{array}{ccc} 
2 & 0 & y \\ 0 & 1 & 0 \\ -\frac{1}{y} & 0 & 0 \end{array}\right):\hs y \neq 0\right\}.
\]
\eex

\subsection{A desingularization of $F_{t, w}$}\lb{subsec-smooth-Xtw}
Recall that if $V$ is an irreducible variety and $V_s$ its set of
smooth points, a {\it desingularization} of $V$ is a proper morphism
$\xi: X \to V$, where $X$ is an irreducible smooth variety and 
$\xi$ maps $\xi^{-1}(V_s)$
isomorphically to $V_s$.

We again assume that $G$ is simply connected. We show that
for any $t \in H$ and $w \in W$, $X_{t, w}$ is a
 desingularization of $F_{t, w}$.

\bpr{pr-Xt-to-Ft} (Grothendieck, see \cite[Theorem 4.4]{slodowy} and 
\cite[Corollary 6.4]{steinberg-inventiones})  For $t \in H$, 
$\mu: X_t \to F_t$ is a {\it desingularization} 
(called the Springer resolution)
of $F_t$. 
\epr

\bth{th-Xtw-smooth} 
For any $t \in H$ and $w \in W$, $X_{t,w}$ is smooth and irreducible,
and $\dim X_{t,w} = \dim G - \dim H - l(w)$. 
\eth

To prove \thref{th-Xtw-smooth}, we need a standard result from algebraic
geometry.

\ble{le-implicitfunction}
Let $f:X\to Z$ be a morphism of smooth algebraic varieties, and
let $Y \subset Z$ be a smooth irreducible subvariety. Suppose 
$f(T_x(X)) + T_{f(x)}(Y) = T_{f(x)}(Z)$ for all $x\in f^{-1}(Y)$. 
Then 
$f^{-1}(Y)$ is smooth and each connected component of $f^{-1}(Y)$
has dimension $\dim X + \dim Y - \dim Z$.
\ele

\noindent
{\it Proof of \thref{th-Xtw-smooth}.}
Consider 
$\mu_t :=\mu|_{X_t}: X_t \to G$. Let $x \in X_{t,w}$, let
$y = \mu(x)$. and let $C_y$ be
the conjugacy class of $y$ in $G$. Then $\mu_t(T_x X_t) \supset T_yC_y$. Since
$C_y$ and $BwB_-$ intersect transversally at $y$, one has
\[
\mu_t(T_x(X_t)) + T_y(BwB_-)\supset T_yC_y + T_y(BwB_-)=T_y(G).
\]
Thus $\mu_t(T_x(X_t)) + T_y(BwB_-)=T_yG$. It follows from \leref{le-implicitfunction} 
that $X_{t,w}$ is smooth with every irreducible component having dimension  
$\dim X_t - l(w)$. 
It remains to show that $X_{t, w}$ is 
irreducible. Let $R_t$ be the regular conjugacy class in $F_t$. 
Since $\mu_t: \mu_t^{-1}(R_t)\to R_t$ is bijective, 
$\mu_t^{-1}(R_t) \cap X_{t,w}$ is irreducible by \prref{pr-smoothirredint}.
For any conjugacy class $C \subset F_t \backslash R_t$, since
$\dim \mu_t^{-1}(C) < \dim X_t$, by dividing $\mu_t^{-1}(C)$ into finitely many 
locally closed $G$-invariant smooth subvarieties of $X_t$, it follows from
\leref{le-implicitfunction}   that 
$\dim (\mu_t^{-1}(C) \cap \mu^{-1}(BwB_-)) < \dim X_t - l(w)$.
By \leref{le-irred},
$X_{t,w}$ is irreducible. Since $\dim X_t = \dim G
- \dim H$, the dimension assertion is clear.
\begin{flushright}
$\square$
\end{flushright}

\bco{co-Xtw-Ftw-resol} For any $t \in H$ and $w \in W$,
$\mu: X_{t, w} \to F_{t, w}$ is a desingularization.
\eco

\begin{proof} This follows directly from \thref{th-Xtw-smooth}
and the definition of a desingularization.
\end{proof}

\subsection{Irreducibility of $X_{t,w}$ for general $G$}

In this subsection, we assume $G$ is connected 
and semisimple but not necessarily
simply connected.

Let $p:\tG \to G$ be a simply connected covering of $G$
 with kernel $Z$. Let $\tB \supset \tT$ be a Borel
subgroup and maximal torus of $\tG$ and assume $p(\tB)=B$ and
$p(\tT)=H$. Let $\tBminus$ be the opposite Borel subgroup such
that $p(\tBminus)=B_-$, and note that $p:\tU \to U$ is an isomorphism,
where $\tU$ is the unipotent radical of $\tB$. For $\tt \in \tT$,
let $\tXtt = \tG \times_{\tB} \tt \tU $ with $\tmu : \tXtt \to \tG$
defined as in \eqref{eq-mu}. For $w\in W$, let
$\tXttw = \tXtt \cap {\tmu}^{-1}(\tB w \tBminus)$.
For $t = p(\tt) \in H$, define 
\[
\tP : \tXtt \to X_t: \; \; \tP [\tg, \tt \tu ]=[p(\tg),p(\tt)p(\tu)].
\]

\ble{le-tXtwirreducible} For any $\tt \in \tT$, $t=p(\tt) \in H$, and
$w\in W$, 

1) $\tP : \tXtt \to X_t$ is an isomorphism of varieties;

2) $\tP (\tXttw)=X_{t,w}$.
\ele

\begin{proof}
For $[p(\tg),p(\tt \tu)] \in G \times_B B$,
${\tP}^{-1}([p(\tg),p(\tt \tu)]) = \{ [\tg, \tt z \tu ]: z\in Z \}$.
Since $[\tg, \tt z \tu] \in \tXtt$ if and only
if $z$ is the identity, $\tP$ is a bijection, so 1) follows since
$X_t$ is smooth. Part 2) is a consequence
of the easily verified assertion
\[
\tmu ([\tg, \tt \tu ]) \in \tB w \tBminus \iff 
\mu(\tP [\tg, \tt \tu ]) \in BwB_-.
\]
\end{proof}

\leref{le-tXtwirreducible} and \thref{th-Xtw-smooth} immediately give

\bpr{pr-Xtw-smooth-any-G}
For a connected complex semisimple Lie group $G$ and for any $t \in H$ and $w \in W$,
$X_{t,w}$ is smooth and irreducible and $\dim X_{t, w} =\dim G - \dim H - l(w)$.
\epr

\sectionnew{The Poisson structure $\Pi$ on $G \times_B B$}\lb{sec-pi}
In this section, $G$ is assumed to be connected and semisimple but not necessarily
simply connected. We will construct and study a Poisson structure 
$\Pi$ on $G \times_B B$ with 
the following properties:

1) The Grothendick map 
$\mu: (G \times_B B, \; \Pi) \to (G, \pi)$ is a Poisson map;

2) $\Pi$ is $H$-invariant, where $H$ acts on $G \times_B B$ by
\eqref{eq-sigma}.

\noindent
Moreover, the $H$-orbits of symplectic leaves of
$\Pi$ in $G \times_B B$ are  
precisely the $X_{t, w}$'s for $t \in H$ and $w \in W$, as defined in
\notaref{nota-Ftw-Xtw}.

\subsection{Definition of the Poisson structure $\Pi$ on $G \times_B B$}\lb{subsec-de-pi}

Recall that $\pi_{\scriptscriptstyle D}^{+}$ is the Poisson structure on $D=G \times G$ 
given in \eqref{eq-piD-plus}. Let $\Bdia = \{(b, b): b \in B\}$, and let
$\phi$ be the projection
\[
\phi: \; \; G \times G \lra (G \times G)/\Bdia.
\]
Recall that $\piG$ is the multiplicative Poisson structure on $G$ defined in 
$\S$\ref{sec-piG}.

\bpr{le-phipi}
$\phi(\pi_{\scriptscriptstyle D}^{+})$ is a well-defined Poisson structure 
on $(G \times G)/\Bdia$, and 
the action of $(\Gdia, \; \piG)$ on 
$((G \times G)/\Bdia, \; \phi(\pi_{\scriptscriptstyle D}^{+}))$ by
left multiplication is   Poisson.
\epr

\begin{proof}
It is well-known \cite{goodearl-yakimov, k-s:quantum} that $\Bdia$ is a Poisson subgroup
of $(\Gdia, \, \piG)$. By \leref{le-D-pipm} in the Appendix, the action of 
$(\Bdia, -\piG)$ on $(G \times G, \pi_{\scriptscriptstyle D}^{+})$
by right multiplication is Poisson. Thus  
$\phi(\pi_{\scriptscriptstyle D}^{+})$ is a well-defined
Poisson structure on $(G \times G)/\Bdia$ by \cite{lu-we:poi}.
 Again by \leref{le-D-pipm} in the Appendix, 
the action by left multiplication of $(\Gdia,  \piG)$ on $((G \times G)/\Bdia, \; 
\phi(\pi_{\scriptscriptstyle D}^{+}))$ is  Poisson.
\end{proof}

Consider the embedding
\begin{equation}\lb{eq-psi-1}
\psi: \; \; G \times_B B \lra (G \times G)/\Bdia:
 \; \; \; [g, \, b] \longmapsto (gb, \, g)_\cdot \Bdia, \; \; g \in G, b \in B.
\end{equation}
Note that $\psi(G \times_B B) = Q/\Bdia$, where
\[
Q  = \{(gb, \,g): g \in G, b \in B\} \subset G \times G.
\]
We show that $Q/\Bdia$ is a Poisson submanifold of $((G \times G)/\Bdia, \, 
\phi(\pi_{\scriptscriptstyle D}^{+}))$ and 
thus obtain a Poisson structure on $G \times_B B$.

For $t\in H$, let $Q_t =\{(gtn, \, g): \; g \in G, n \in U\}.$
Then $Q_t$ is $\Bdia$-invariant, and
$\psi(X_t) =Q_t/\Bdia$.

\bpr{pr-Qt-sub} For every $t \in H$, $Q_t/\Bdia$ is a Poisson submanifold of 
$(G \times G)/\Bdia$ with respect 
to $\phi(\pi_{\scriptscriptstyle D}^{+})$.  
\epr

\begin{proof}  
By \prref{pr-appendix-sub} in the Appendix, it suffices to show that
$Q_t$ is coisotropic in $G \times G$ with respect to 
$\pi_{\scriptscriptstyle D}^{+}$ and that the characteristic distribution
(see below) 
of $\pi_{\scriptscriptstyle D}^{+}$ 
on $Q_t$ coincides with the distribution defined by the tangent spaces to the
$\Bdia$-orbits in $Q_t$. 

By the definition of $\pi$ in \notaref{nota-pi0}, 
\[
\eta: \; \; \; (G \times G, \, \pi_{\scriptscriptstyle D}^{+})\lra (G, \, \pi):   \;\;\;
(g_1, g_2) \longmapsto g_{1} g_2^{-1}
\]
is Poisson. Let $\tau: G \times G \to G \times G: (g_1, g_2)\mapsto (g_1^{-1}, g_2^{-1})$. 
It is clear 
from the definition of
$\pi_{\scriptscriptstyle D}^{+}$ that $\tau$ preserves $\pi_{\scriptscriptstyle D}^{+}$. 
Thus 
\[
\eta_1 := \eta  \tau: \; \; \; (G \times G, \, 
\pi_{\scriptscriptstyle D}^{+})\lra (G, \, \pi): 
\;\;\;   (g_1, g_2) \longmapsto g_{1}^{-1} g_2
\]
is Poisson. 
Let $t \in H$.
By \leref{le-tN-coisotropic}, $Ut^{-1}= t^{-1}U$ is
coisotropic in $G$ with respect to  $\pi$. 
Since 
$\eta_1$ is a  
submersion, it follows from
\cite{we:coiso} that $Q_t=\eta_1^{-1}(Ut^{-1})$ is coisotropic in $D$ with respect to 
$\pi_{\scriptscriptstyle D}^{+}$.
 
Recall from $\S$\ref{subsec-coiso-reduction}
in the Appendix that the characteristic distribution of
$\pi_{\scriptscriptstyle D}^+$ on 
$Q_t$ is by definition the image of the bundle map 
\[
\tilde{\pi}_{\scriptscriptstyle D}^{+}: \; \; (TQ_t)^0 \lra TQ_t,
\]
where  $(TQ_t)^0$ is the 
conormal bundle of $Q_t $ in $G \times G$, and 
$\tilde{\pi}_{\scriptscriptstyle D}^{+}: T^*(G \times G) \to T(G \times G)$ is
defined as in \eqref{eq-tilde-pi}.
Fix $g \in G$, $n \in U$,  and let $d = (gtn, \, g) \in Q_t$. 
We now use formula \eqref{eq-tilde-piD-1} in the proof of 
\leref{le-leaves-piD} in the Appendix to compute 
$\tilde{\pi}_{{\scriptscriptstyle D}}^{+} ((T_dQ_t)^0)$.
First,
\[
T_dQ_t  = r_d \gdia + l_d (\n \oplus 0) =   
{r_d} \left( \Ad_{(g, g)}(\gdia + (\n \oplus 0))\right).
\]
Identify $\d$ with $\d^*$ via the bilinear form $\lara$, and for $(x, y) \in \d$, let
$\alpha_{(x, y)}$ be the right invariant $1$-form on $G \times G$ 
whose value at the identity
element is $(x, y)$. Then 
\[
(T_dQ_t)^0 = \{\alpha_{(x, x)}(d): \; \; x \in \Ad_g \b\}.
\]
Let $p_1: \d \to \gdia$ be the projection with respect to the decomposition $\d = 
\gdia + \gst$. Let $x \in  \Ad_g \b$. By formula \eqref{eq-tilde-piD-1} in the proof of 
\leref{le-leaves-piD} in the Appendix,
\[
\tilde{\pi}_{{\scriptscriptstyle D}}^{+}(\alpha_{(x, x)}(d)) = l_d p_1 \Ad_{d}^{-1}(x, x)
=l_d p_1 (\Ad_{tn}^{-1}\Ad_g^{-1}x, \; \Ad_g^{-1}x)=l_d(\Ad_g^{-1}x, \; \Ad_g^{-1}x),
\]
which is exactly the infinitesimal generator of the $\Bdia$-action on $Q_t$
in the direction of $(\Ad_g^{-1}x, \; \Ad_g^{-1}x)$. This finishes the proof of
\prref{pr-Qt-sub}.
\end{proof}

\bco{co-Q-sub}
$Q/\Bdia$ is a Poisson submanifold of $(G \times G)/\Bdia$ 
with respect to the Poisson structure 
$\phi(\pi_{\scriptscriptstyle D}^{+})$ in \prref{le-phipi}.
\eco

\begin{proof} \coref{co-Q-sub} follows immediately from \prref{pr-Qt-sub} and
the fact that $Q/\Bdia = \bigcup_{t \in H} (Q_t/\Bdia)$.
\end{proof}

Recall the isomorphism $\psi: G \times_B B \to Q/\Bdia$ in \eqref{eq-psi-1}.

\bnota{de-pi-GBB}
The Poisson structure $\psi^{-1} (\phi(\pi_{\scriptscriptstyle D}^+))$ 
on $G \times_B B$ will 
be denoted by $\Pi$.
\enota
 
We summarize some of the properties of the Poisson structure $\Pi$ on $G \times_B B$. 
Recall that $\sigma$ is the left action of $G$ on $G \times_B B$ given in \eqref{eq-sigma}.

\bpr{pr-properties-pi} 1) For each $t \in H$, $X_t$ is a 
Poisson submanifold of $(G \times_B B, \Pi)$;

2) The Grothendieck map $\mu: (G \times_B B, \, \Pi) \to (G, \pi)$
is Poisson;

3) With the Poisson structure $\Pi$ on $G \times_B B$, $\sigma$ is a Poisson action
by the Poisson Lie group $(G, \piG)$;

4) $\Pi$ is $H$-invariant for the action $\sigma$ of $H$ on $G \times_B B$. 
\epr

\begin{proof} Part 1) follows from the 
definition of $\Pi$ and \prref{pr-Qt-sub}, and 
 2) is a consequence of 
the definitions of $\Pi$  and $\pi$ and 
\coref{co-Q-sub}.
Part 3) 
follows from the definition of $\Pi$
and \prref{le-phipi}, and 4) is then 
clear since $\piG$ vanishes at points in $H$.
\end{proof}

\subsection{Symplectic leaves of $\Pi$ in $G \times_B B$}\lb{subsec-leaves-pi}

We use \prref{pr-appendix-sub} in the Appendix to
determine the symplectic leaves of $\Pi$ in $G \times_B B$.

\ble{le-leaves-piD-GG} The symplectic leaves of
 $\pi_{{\scriptscriptstyle D}}^{+}$ 
are the connected components of 
the nonempty intersections $G^*(h\dw, e)\Gdia \cap \Gdia (e, h^\prime \du)G^*$,
where $h, h^\prime \in H$ and $w, u \in W$. 
\ele

\begin{proof}
By \leref{le-leaves-piD} in the Appendix, 
symplectic leaves of $\pi_{{\scriptscriptstyle D}}^{+}$
are the connected components of nonempty intersection of
$(G^*, \Gdia)$ and $(\Gdia, G^*)$-double cosets in $G \times G$. 
By \leref{le-Gs-orbits-on-G}, $(G^*, \Gdia)$ and $(\Gdia, G^*)$ double cosets in $G \times G$
are respectively of the form $G^*(h\dw, e)\Gdia$ and $ \Gdia (e, h^\prime \du)G^*$,
where $h, h^\prime \in H$ and $w, u \in W$. 
\end{proof}

Recall that for $h \in H$ and $w \in W$, $G^*{}_\cdot (h\dw)$ is the 
$G^*$-orbit in $G$ through $h \dw$, where $G^*$ acts on $G$ by \eqref{Gs-on-G}. 

\bpr{pr-leaves-pi} 1) For any $t, h \in H$ and $w \in W$, 
$X_t \cap \mu^{-1}(G^*{}_\cdot (h\dw))$ is nonempty, smooth, 
and all of its connected components have dimension 
\[
\dim G - \dim H - l(w) - \dim(H/H_w);
\]

2) The symplectic leaves of $\Pi$ in $G \times_B B$ 
are precisely the connected components
of $X_t \cap \mu^{-1}(G^*{}_\cdot (h\dw))$, where $t, h \in H, w \in W$.
\epr

\begin{proof} 1) Let $R_t$ be the unique regular conjugacy class in $F_t$. By \coref{R},
$R_t \cap (G^*{}_\cdot (h\dw)) \neq \emptyset$. Since $\mu|_{X_t}: X_t \to 
F_t$ is onto, $X_t \cap \mu^{-1}(G^*{}_\cdot (h\dw))\neq \emptyset$. 
An argument similar to the  proof of the smoothness of $X_{t,w}$ in
\thref{th-Xtw-smooth}
shows that $X_t \cap \mu^{-1}(G^*{}_\cdot (h\dw))$ is smooth. The dimension
assertion follows from \leref{le-implicitfunction}.

2) Let $t \in H$. Let $t_1 \in H$
be such that $t_1^2 = t$. Then for any $g \in G$ and $n \in U$, 
$(gtn, g) = (gt_1, gt_1)(t_1n, t_1^{-1}) \in \Gdia G^*$. 
Thus $Q_t \subset \Gdia G^*$.
By \prref{pr-appendix-sub} in the Appendix, symplectic leaves of 
$\phi(\pi_{\scriptscriptstyle D}^+)$ in
$Q_t/\Bdia$ are the 
connected components of $(Q_t \cap (G^*(h\dw, e)\Gdia))/\Bdia$, where $w \in W$ and $h \in H$.
It is routine to check that $\mu^{-1}(G^*{}_\cdot (h\dw))
=\psi^{-1}((G^*(h\dw, e)\Gdia))/\Bdia))$.
Now 2) follows
from the isomorphism $\psi: G \times_B B \to Q/\Bdia$  and the 
fact the $X_t$'s are
Poisson submanifolds of $(G \times_B B, \, \Pi)$.
\end{proof}


\subsection{$H$-orbits of symplectic leaves of $\Pi$ in $G \times_B B$}\lb{subsec-H-leaves-pi}

By 4) of \prref{pr-properties-pi}, the Poisson structure $\Pi$ on
$G \times_B B$ is $H$-invariant for the $H$-action in \eqref{eq-sigma}.

\bth{th-H-leaves-Pi}
The $H$-orbits of symplectic leaves of $\Pi$ in $G \times_B B$ are precisely the 
irreducible smooth subvarieties $X_{t, w}$, where $t \in H$ and $w \in W$. The
dimension of every symplectic leaf in $X_{t,w}$ is $\dim G - \dim H - l(w) - \dim(H/H_w)$. 
\eth

\begin{proof}
By \prref{pr-leaves-pi} and the fact that each $X_{t,w}$ is $H$-invariant,
 every $H$-orbit of symplectic leaves of $\Pi$ in $G \times_B B$
is contained in $X_{t, w}$ for some $t \in H$ and $w \in W$. Fix
$t \in H$ and $w \in W$. 
Then $X_{t, w}$ is nonempty, smooth and
connected by \prref{pr-Xtw-smooth-any-G}. 
To prove \thref{th-H-leaves-Pi}, it suffices to show that
for any connected component $S$ in $X_t \cap \mu^{-1}(G^*{}_\cdot (h \dw))$, where
$h \in H$, the map 
\begin{equation}\lb{sigma1}
\sigma_1: \; \; H \times S \longrightarrow X_{t, w}: \; \; \; (h, \, [g, b]) \longmapsto
[hg, b], \hs h \in H, [g, b] \in S
\end{equation}
has surjective differential everywhere. Let
\[
Q_{t, w} = Q_t \cap ((B \times B_-) (h\dw, e) \Gdia ).
\]
Then $Q_{t, w}$ is smooth and irreducible since it is a principal 
$B$-bundle over  the smooth irreducible variety $X_{t, w}$.
 An argument similar to an argument in 
the proof of  \prref{pr-H-orbits}, Part 2),
 shows that for any open submanifold
$S^\prime$ of $Q_t \cap G^* (h \dw, e) \Gdia$, the map
\[
H \times S^\prime \longrightarrow  Q_{t,w}: \; \; 
(h, (g_1, g_2)) \longmapsto (hg_1, hg_2), \hs h \in H, (g_1, g_2) \in S^\prime
\]
has surjective differential everywhere. It follows that $\alpha_1$ in
\eqref{sigma1} has surjective differential everywhere.
\end{proof}

\subsection{$H$-equivariant Poisson desingularization}\lb{poissondesing}

An irreducible Poisson variety $(X,\pi)$ is said to be regular if 
$\pi$ has constant rank on $X$.

\bde{de-poissondesingularization} 
Let $(X,\piX)$ and $(Y,\piY)$ be irreducible Poisson varieties, with $Y$
normal. A desingularization $f: X \to Y$ is called a  {\it Poisson 
desingularization} if 
$\piX$ is a regular Poisson structure on $X$ and if $f: (X, \piX) \to (Y, \piY)$ is Poisson.
If, in addition, a torus $H$ acts by Poisson isomorphisms on 
both $(X, \piX)$ and $(Y, \piY)$ such that $f$ is $H$-equivariant and 
$H$ acts transitively on the set of symplectic leaves of $\piX$ in $X$, we
call $f: (X, \piX) \to (Y, \piY)$ an {\it $H$-equivariant Poisson
desingularization}.
\ede

\bre{re-symplecticresolution}
If $\pi_{\scriptscriptstyle Y}$ is generically nondegenerate,  
$f:X \to Y$ is a symplectic resolution in the
sense of \cite{kasym} (see also \cite{beauville,fu-symplectic,fu-survey}). 
Our notion of Poisson desingularization
is different from that of Poisson resolution  by Fu in \cite{fu-poisson}. 
It would be interesting to relate ideas from this paper to 
the recent paper \cite{camille} by Laurent-Gengoux.
\ere

\bco{co-poissondesingularization}
For any $t \in H$ and $w \in W$, 
$\mu: (X_{t,w}, \Pi) \to (F_{t, w}, \pi)$ is 
an $H$-equivariant Poisson desingularization.
\eco

\bex{ex-SL2-Xtw-Ftw} For $G = SL(2, \Cset)$, $X = G \times_B B$ 
is the union of two open charts $X^\prime$ and $X^{\prime\prime}$, where
\begin{align*}
X^\prime &= \left\{\left[ \left(\begin{array}{cc} x_1 & -1 \\ 1 & 0 \end{array}
\right), \, \left(\begin{array}{cc} h_1 & y_1 \\ 0 & \frac{1}{h_1} \end{array}\right)\right]:
\; x_1, y_1, h_1 \in \Cset, h_1 \neq 0\right\}\\
X^{\prime\prime} &= \left\{\left[ \left(\begin{array}{cc} 1 & 0 \\ x_2 & 1\end{array}
\right), \, \left(\begin{array}{cc} h_2 & y_2 \\ 0 & \frac{1}{h_2} \end{array}\right)\right]:
\; x_2, y_2, h_2 \in \Cset, h_2 \neq 0\right\},
\end{align*}
and 
$h_1 = h_2,  x_2 = \frac{1}{x_1}$ and $y_2 = x_1 (h_1 - \frac{1}{h_1} + x_1 y_1)$
on $X^\prime \cap X^{\prime\prime}$.
Using the formula for $\pi$ on $G$ in \exref{ex-sl2-pi0} and the
fact that $\mu: (X, \Pi) \to (G, \pi)$ is Poisson, one can compute the
Poisson structure $\Pi$ on $X$ and get
\begin{align}\lb{eq-Pi-1}
\{h_1, x_1\}_\Pi = \{h_1, y_1\}_\Pi = 0, \hs \{x_1, y_1\}_\Pi
= h_1 + x_1 y_1 \hs \mbox{on} \; \; X^\prime, \hs \mbox{or}\\
\lb{eq-Pi-2}
\{h_2, x_2\}_\Pi = \{h_2, y_2\}_\Pi = 0, \hs \{x_2, y_2\}_\Pi
= -\left(\frac{1}{h_2} + x_2 y_2\right) \hs \mbox{on} \; \; X^{\prime \prime}.
\end{align}
It follows that each $X_t$ is a Poisson variety of $(X, \Pi)$. Moreover, 
for  $t = {\rm diag}(h, \frac{1}{h})$, note that $X_{t, w_0} 
\subset X^\prime \cap X^{\prime\prime}$ and is given by $h + x_1 y_1 =0$ or
$\frac{1}{h} + x_2 y_2 = 0$,
so it is clear from \eqref{eq-Pi-1} and \eqref{eq-Pi-2} that $\Pi$ vanishes on $X_{t, w_0}$
and is symplectic on $X_{t, w=1}$. On the other hand, $\mu: X_{t, w_0} \to
F_{t, w_0}$ is an isomorphism for any $t$, and
so is
$\mu: X_{t, w=1} \to F_{t, w=1}$ unless $h = \pm 1$, in which case $F_{t, w=1}$ 
is singular at $t$ and $\pi$ has two $H$-orbits of symplectic leaves in $F_{t, w=1}$.
\eex

\subsection{Remarks on the Poisson structure
$\phi(\pi_{\scriptscriptstyle D}^+)$ on $(G \times G)/\Bdia$}\lb{remarksD}

Consider the isomorphism 
\begin{equation}\lb{eq-gamma}
\gamma: \; \; (G \times G)/\Bdia \lra (G/B) \times G: \; \; (g_1, g_2)_\cdot \Bdia 
\longmapsto (g_1 {}_\cdot B, \, g_1g_2^{-1}).
\end{equation}
Composing $\gamma$ with the embedding $\psi: G \times_B B \to (G \times G)/\Bdia$ in 
\eqref{eq-psi-1} gives 
\[
\gamma \psi: \; G \times_B B \lra (G/B) \times G: \; \; [g, b] \longmapsto 
(g_\cdot B, \, gbg^{-1}), \; \; g \in G, b \in B
\]
with 
\[
\gamma \psi (G \times_B B) = \{(g_1 {}_\cdot B, g_2): g_1, g_2 \in G, g_1^{-1} g_2 g_1 \in B\}.
\]
Using the embedding $\gamma \psi$, we can regard $(G \times_B B, \Pi)$ as a Poisson
submanifold of $((G/B) \times G, \,\gamma \phi(\pi_{\scriptscriptstyle D}^+))$,
where recall from $\S$\ref{subsec-de-pi} 
that $\phi$ is the projection $G \times G \to (G \times G)/\Bdia$.
 
In this subsection, we prove some properties of the 
Poisson structure $\gamma \phi(\pi_{\scriptscriptstyle D}^+)$ on $(G/B) \times G$
that are helpful in understanding the Poisson structure $\Pi$ on $G \times_B B$.
For notational simplicity, set
\[
\Pi_1 = \gamma \phi(\pi_{\scriptscriptstyle D}^+).
\]

Let $\zeta: G \to G/B$ be the projection.
Since $B$ is a Poisson subgroup of $(G, \piG)$, 
$\zeta(\piG)$ is a well-defined Poisson structure on $G/B$, which
we denote by $\pi_{\scriptscriptstyle G/B}$,
and $(G/B, \pi_{\scriptscriptstyle G/B})$ is 
a Poisson $(G, \piG)$-homogeneous space. 
In particular, $\pi_{\scriptscriptstyle G/B}$ is invariant under
translation by elements in $H$. 
It is shown in \cite{goodearl-yakimov} that
the $H$-orbits of symplectic leaves of $\pi_{\scriptscriptstyle G/B}$ 
in $G/B$ are precisely the intersections
$(Bu_\cdot B) \cap (B_- v_\cdot B)$ for $u, v \in W, v \leq u$. 
In particular, $\pi_{\scriptscriptstyle G/B}$ vanishes at $u_\cdot B$ for every $u \in W$.

\bpr{p1p2}
Both of the following projections are Poisson:
\begin{align*}
p_1:& \; \; ((G/B) \times G, \, \,\Pi_1) 
\lra (G/B, \, \pi_{\scriptscriptstyle G/B}): \; \; (g_1{}_\cdot B, \; g_2)
 \longmapsto
g_1{}_\cdot B\\
p_2: & \; \; ((G/B) \times G, \,\, \Pi_1) \lra (G, \, \pi): \; \; (g_1{}_\cdot B, \;g_2) 
\longmapsto g_2.
\end{align*}
\epr

\begin{proof} The map $p_2 \gamma \phi: G \times G \to G$ is $\eta$ as in \eqref{eq-eta},
so $p_2(\Pi_1) = \eta(\pi_{\scriptscriptstyle D}^+)=\pi$ by the definition of $\pi$.
Let $p_1^\prime: G \times G \to G$ be the projection to the first factor. Then
$p_1(\Pi_1) = (p_1 \gamma \phi)(\pi_{\scriptscriptstyle D}^+) = 
(\zeta p_1^\prime)(\pi_{\scriptscriptstyle D}^+)$. Note that $p_1^\prime$ is a group
homomorphism, and for $\Lambda \in \wedge^2(\g \oplus \g)$ in 
 \eqref{eq-Lambda},
$p_1^\prime (\Lambda) = \Lambda_0$, where $\Lambda_0 \in \wedge^2 \g$ is given in 
\prref{pr-piG}. Thus $p_1^\prime(\pi_{\scriptscriptstyle D}^+) = 
\Lambda_0^R + \Lambda_0^L$ (see $\S$\ref{intro-notation}), and
\[
p_1(\Pi_1) = (\zeta p_1^\prime)(\pi_{\scriptscriptstyle D}^+) 
= \zeta(\Lambda_0^R + \Lambda_0^L)
=\zeta(\Lambda_0^R) = \pi_{\scriptscriptstyle G/B}.
\]
\end{proof}

\bco{co-piGB} The projection $q: (G  \times_B B, \, \Pi) 
\to (G/B, \, \pi_{\scriptscriptstyle G/B}): \, [g, b] \mapsto g_\cdot B$
is Poisson.
\eco

\sectionnew{Birational isomorphisms between $X_{t, w}, F_{t, w}$
and $G^{1, w^{-1}w_0} \times_H G^{1, w_0}/H$}\lb{sec-log-pi}


\subsection{Relation between $X_{t,w_0}$ and $X_{t, w}$}\lb{subsec-tw-tw0}

Let $w_0$ be
the longest element in $W$. 
Recall that $\sigma: G \times (G \times_B B) \to G \times_B B$ is the action 
of $G$ on $G \times_B B$ given by \eqref{eq-sigma} and that
$G^{u, v} = BuB \cap B_- v B_- \subset G$ for $u, v \in W$.

\ble{le-key}
For any $w \in W$,
\[
\sigma\left(G^{1, w^{-1}w_0} \, \times \, X_{t, w_0} \right) \subset X_{t,w}.
\]
\ele

\begin{proof} 
Let $[g,tn] \in X_{t, w_0}$, so $g\in G$, $n\in U$, and 
$gtng^{-1} \in Bw_0B_- =B w_0$. 
If $g_1 \in G^{1, w^{-1}w_0}$, then
 $g_1^{-1} \in B_- w_0^{-1} w B_-$. Thus 
\[
\mu ([g_1g, b])= g_1 gtng^{-1} g_1^{-1} \in Bw_0 B_- w_0^{-1} w B_- = BwB_-.
\]
\end{proof} 

By \leref{le-leaves-piG}  and \thref{th-H-leaves-Pi},
$G^{1, w^{-1}w_0}$ and $X_{t,w_0}$ are Poisson submanifolds of
$(G,\piG)$ and $(X_t,\Pi)$ respectively.
\prref{pr-properties-pi} implies the following

\bpr{pr-key}
For $w \in W$, equip $G^{1, w^{-1}w_0}$ with the Poisson structure $\piG$ and $X_{t, w}$ the
Poisson structure $\Pi$. Then
\[
\sigma: \; \; (G^{1, w^{-1}w_0}, \piG) \, \times \, (X_{t, w_0}, \Pi)
 \lra  (X_{t,w}, \Pi)
\]
is a Poisson map.
\epr

\subsection{A decomposition of $X_{t, w_0}$}\lb{subsec-Xtw0}
We partition
$X_{t, w_0}$ into smooth subvarieties.
Recall that $q: G \times_B B \to G/B: [g, b] \mapsto g_\cdot B$ is the projection. 
Fix $t \in H$. Set
\begin{equation}\lb{eq-Xutw0}
X^u_{t, w_0} = X_{t, w_0} \cap q^{-1}(B w_0 u \cdot B) \hs \mbox{for} \; \; u \in W.
\end{equation}
The Bruhat decomposition gives
\begin{equation}\lb{Xtw0-Xutw0}
X_{t, w_0} = \bigsqcup_{u \in W} X^u_{t, w_0} \hs (\mbox{disjoint union}).
\end{equation}

\ble{le-Xutw0-empty}
Let $t \in H$ and $u \in W$. If $X^u_{t, w_0} \neq \emptyset$, then $u \leq w_0 u$,
where $\leq$ is the Bruhat order on $W$, and
\[
q\left(\Xu\right) \subset (B_-u_\cdot B) \cap (B w_0 u_\cdot B).
\]
\ele

\begin{proof} 
Assume that $[g, tn] \in \Xu$, where $g \in G$ and $n \in U$. 
Since $g \in B w_0uB$ and
$g tng^{-1} \in Bw_0$, 
$g \in (Bw_0)^{-1} B w_0 u B = B_- u B$.
Thus $g \in B_- u B \cap Bw_0uB$
and $u \leq w_0u$ by \cite[Corollary 1.2]{deo}.
\end{proof}

For $u \in W$ such that $u \leq w_0 u$, set
\[
{\calR}_{u, w_0u} = (B_-u_\cdot B) \cap (B w_0 u_\cdot B) \subset G/B.
\]
We now establish an isomorphism between 
$X^u_{t, w_0}$ and ${\calR}_{u, w_0u} \times U \cap u^{-1}U_- u$.
To this end, first parametrize $Bw_0u_\cdot B$ by the isomorphism
\[
U_{w_0u} \lra B w_0 u_\cdot B:\; \; \;  m \longmapsto m \dw_0 \du_\cdot B,
\]
where $U_{w_0u} = U  \cap (w_0uU_- u^{-1}w_0^{-1})$. Note that
for $m 
\in U_{w_0u}$, $m \dw_0 \du_\cdot B \in B_-u_\cdot B$ if and only if
$m\in U_{w_0u} \cap B_-uB u^{-1} w_{0}^{-1}$. Thus
\begin{equation}\lb{R-iso}
U_{w_0u} \cap B_-uB u^{-1} w_{0}^{-1} \lra {\calR}_{u, w_0u}: \; \; 
m \longmapsto m\dw_0 \du_\cdot B
\end{equation}
is an isomorphism. Let $m\in U_{w_0u} \cap B_-uB u^{-1} w_{0}^{-1}$. By
the unique factorization 
\[
B_- \times (U \cap u^{-1}U u) \lra B_- u B, \, \, (b_-, n)\longmapsto b_-\du n, \; \; \;
b_- \in B_-, \, n \in U \cap u^{-1}U u,
\]
\begin{equation}\lb{nm}
m\dw_0 \du = b_-\du n_m, \hs  \mbox{for unique}\; \;b_-\in B_- \; \mbox{and}\;
n_m \in U \cap u^{-1}U u.
\end{equation}
Define
\[
\xi^u: \; \; U_{w_0u} \cap B_-uB u^{-1} w_{0}^{-1} 
\lra U \cap u^{-1}U u: \; \; m \longmapsto n_m.
\]

\bnota{nota-Uv-G0}
Let $G_0 = B_-B$. For $g \in B_- B$, write
\[
g = (g)_- (g)_0 (g)_+, \hs \mbox{where}\; (g)_- \in U_-, (g)_0 \in H, (g)_+ \in U.
\]
\enota

For $m \in U_{w_0u} \cap B_-uB u^{-1} w_{0}^{-1}$ with the decomposition 
of $m\dw_0 \du$ in \eqref{nm}, since
$m\dw_0 = b_- \du n_m \du^{-1} \in G_0$, $n_m = \du^{-1} (m\dw_0)_+ \du$. Thus the
map $\xi^u$ is also given by
\begin{equation}\lb{eq-xiu}
\xi^u: \; \; U_{w_0u} \cap B_-uB u^{-1} w_{0}^{-1} 
\lra U \cap u^{-1}U u: \; \; m \longmapsto
\du^{-1} (m\dw_0)_+ \du.
\end{equation}
Consider now the morphism
\[
J^u_t: \; {\calR}_{u, w_0u} \times (U \cap u^{-1}U_- u) \lra X_{t}: \; 
(m \dw_0 \du_\cdot B, \, m_1) \longmapsto [m \dw_0 \du, \,\, t m_1 \xi^u(m)]
\]
where $m \in U_{w_0u} \cap B_-uB u^{-1} w_{0}^{-1}$ and $m_1 \in U \cap u^{-1}U_- u.$

\ble{le-Jut} The image of $J^u_t$ is contained in $X_{t, w_0}^u$.
\ele
 
\begin{proof}
Let $m \in U_{w_0u} \cap B_-uB u^{-1} w_{0}^{-1}$, and write $m \dw_0 \du = b_- \du \xi^u(m)$ for
unique $b_- \in B_-$. Then
for any $m_1 \in U \cap \du^{-1}U_- \du$,
\[
(m \dw_0 \du) (tm_1 \xi^u(m)) (m \dw_0 \du)^{-1} =  m \dw_0 \du tm_1 \xi^u(m) \xi^u(m)^{-1} 
\du^{-1} b_-^{-1} \in Bw_0B_- = Bw_0.
\]
Thus $[m \dw_0 \du, \,\, t m_1 \xi^u(m)] \in X_{t, w_0} \cap q^{-1}(Bw_0u_\cdot B)=
X_{t, w_0}^u$. 
\end{proof}

\bpr{pr-Xu}
For any $u \in W$ such that $u \leq w_0u$, 
\begin{equation}\lb{eq-Ju}
J^u_t: \; {\calR}_{u, w_0u} \times (U \cap u^{-1}U_- u) \lra X_{t,w_0}^u
\end{equation}
is an isomorphism. In particular, $\Xu$ is smooth and irreducible and 
\[
\dim \Xu = l(w_0) - l(u).
\]
\epr

\begin{proof} 
Since $U_{w_0u} \times U \to X_t:   (m, \, n) \longmapsto [m \dw_0 \du, \; tn]$
is an embedding, $J^u_t$ is an embedding. Define 
$K_t^u: \Xu \to {\calR}_{u, w_0u} \times (U \cap u^{-1}U_- u)$ as follows:
Let $[g,tn] \in \Xu$. 
Assume without loss of generality that $g = m \dw_0 \du$ with $m  \in 
U_{w_0u} \cap B_-uB u^{-1} w_{0}^{-1}$ and write 
$n = m_1 n_1$ with $m_1 \in U \cap u^{-1}U_- u$ and $n_1 \in U \cap u^{-1}U u.$
Let
\[
K_t^u ([m \dw_0 \du, tm_1n_1]) = (m \dw_0 \du_\cdot B, m_1).
\]
It is straightforward to check that $K_t^u$ and $J_t^u$ are inverse isomorphisms.
By \cite{KL}, ${\calR}_{u, w_0u}$ is smooth  and irreducible and has dimension
$l(w_0)-2l(u)$.  
It follows that 
$\Xu$ is smooth and irreducible, and that
$\dim \Xu = l(w_0)-l(u).$
\end{proof}

\subsection{The open subvariety $X_{t,w_0}^{1}$ of $X_{t, w_0}$}\lb{subsec-Xw01}
For $u = 1 \in W$, the open subset $X_{t,w_0}^{1}$ of $X_{t, w_0}$
is especially simple.
Indeed, $R_{1, w_0}= {B_-}_\cdot B \cap {Bw_0}_\cdot B$ is parametrized by
\begin{equation}\lb{eq-Rw0}
U \cap B_- w_0 B_- \lra {\calR}_{1, w_0}:\;\;\;  
m \longmapsto {m\dw_0}_\cdot B,
\end{equation}
and the isomorphism $J_t := J_t^1$ in \eqref{eq-Ju} simplifies to 
\begin{equation}\lb{eq-Jt}
J_t: \; \; {\calR}_{1, w_0} \lra X_{t, w_0}^{1}: \, \, 
{m \dw_0}_\cdot B \longmapsto [m \dw_0, \; t (m\dw_0)_+],
\hs m \in  U \cap B_- w_0 B_-
\end{equation}
The inverse of $J_t$ is the restriction to $X_{t, w}^1$ of the projection 
$q: G \times_B B \to G/B$.
Identify $U \cap B_- w_0 B_-$ with $G^{1, w_0}/H$ by $m \mapsto m_\cdot H$
for $m \in U \cap B_- w_0 B_-$. Then the isomorphism in \eqref{eq-Rw0} can be
replaced by
\begin{equation}\lb{eq-psi-w0}
\psi_{w_0}: \; \; G^{1, w_0}/H \lra {\calR}_{1, w_0}:\;\;\;  
g_\cdot H \longmapsto {g\dw_0}_\cdot B.
\end{equation}

\ble{le-comp} The composition $J_t \psi_{w_0}$ is given by
\[
J_t \psi_{w_0}: \; \; \; G^{1, w_0}/H \lra X_{t, w_0}^1: \; \; 
g_\cdot H \mapsto [g\dw_0, \, t (g\dw_0)_+].
\]
\ele

\begin{proof} For $g \in G^{1, w_0}$, write $g = mh$, where $m \in U \cap 
(B_- w_0 B_-)$ and $h \in H$. It is straightforward to check that
$[m\dw_0, \, t(m\dw_0)_+] = [g\dw_0, \, t (g\dw_0)_+]$.
\end{proof}

Since $X_{t, w_0}^1$ is open in 
$X_{t, w_0}$, $X_{t, w_0}^1$ is a Poisson subvariety of
$(G \times_B B, \Pi)$.
Recall that $\pi_{\scriptscriptstyle G/B}$
is the projection to $G/B$ of the Poisson structure $\piG$ on $G$, and that
${\calR}_{1, w_0}$ is a Poisson 
subvariety of $(G/B, \piGB)$ by \cite{goodearl-yakimov}. Moreover,  
since the Poisson structure $\piG$ on $G$ is invariant under
right multiplication by elements in $H$, 
the quotient $G/H$ has a well-defined
Poisson structure which we   denote by $\piGH$. Clearly
$G^{u, v}/H$ is a Poisson subvariety of $(G/H, \piGH)$ for any $u, v \in W$.
 
\ble{le-Jt-poi} Both 
\[
\psi_{w_0}: \;  (G^{1, w_0}/H, \,\piGH)  {\lra} ({\calR}_{1, w_0}, \,\piGB)
\hspace{.1in} \mbox{and} \hspace{.1in}
J_t: \;  ({\calR}_{1, w_0}, \,\piGB)  {\lra} (X_{t, w_0}^{1}, \,\Pi)
\]
are Poisson isomorphisms.
\ele

\begin{proof} By $\S$\ref{remarksD},
$\pi_{{\scriptscriptstyle G/B}}$ vanishes at ${ \,\dot{w}_0}_\cdot B \in G/B$.
Since the action of $(G, \piG)$ on $(G/B, \, \piGB)$ by left translation
is Poisson, the map 
\[
(G^{1, w_0}, \; \piG) \lra ({\calR}_{1, w_0}, \; \piGB): \; \; g \longmapsto {g\dw_0}_\cdot B
\]
is Poisson, and so is $\psi_{w_0}$. For any $t \in H$, the projection
$q: (G \times_B tU,,\Pi) \to (G/B, \piGB)$ is Poisson by \coref{co-piGB}. Thus
$J_t^{-1} = q|_{X_{t, w_0}^{1}}:
(X_{t, w_0}^{1}, \, \Pi) \to ({\calR}_{1, w_0}, \piGB)$ is Poisson.
\end{proof}

We now state a fact that will be used in $\S$\ref{subsec-main}. 
Let
\begin{equation}\lb{eq-xi-0}
\xi := \xi^1:\; \;  U \cap B_- w_0 B_- \lra  U: \; \; m \longmapsto (m \dw_0)_+.
\end{equation}
The following \leref{le-xi} can be checked directly.

\ble{le-xi}
The image of $\xi$ in \eqref{eq-xi-0} is again $U \cap B_- w_0 B_-$ and 
\begin{equation}\lb{eq-xi}
\xi: \; \; U \cap B_- w_0 B_- \lra U \cap B_- w_0 B_-: \, \, m \longmapsto (m \dw_0)_+
\end{equation}
is biregular with inverse given by
\begin{equation}\lb{eq-xi-inverse}
\xi^{-1}: \; \; U \cap B_- w_0 B_- \lra U \cap B_- w_0 B_-: \, \, 
n \longmapsto (n \dw_0^{-1})_+.
\end{equation}
\ele

\subsection{A birational isomorphism between $X_{t,w}$ and 
$G^{1, w^{-1}w_0} \times_H G^{1, w_0}/H $}\lb{subsec-main}

For $t \in H$, consider the Zariski open subset $X_t^0$ of $X_t$ given by
\[
X_t^0 = \{[n_1\dw_0, \, tn_2]: \, n_1 \in U, 
n_2 \in U \cap (B_- w_0 B_-)\}.
\]
For $w \in W$, $X_{t, w} \cap X_t^0$ is then a Zariski open
subset of $X_{t, w}$, and we show below that it is nonempty.
In addition, consider the free right action of $H \times H$ on 
$G^{1, w^{-1}w_0} \times G^{1, w_0} $ given by
\begin{equation}\lb{eq-HH-action}
(g_1, g_2 ) \cdot (h_1, h_2) = 
(g_1 h_1, \,\, h_1^{-1} g_2 h_2), \hs g_1, g_2 \in G, h_1, h_2 \in H.
\end{equation}
The action preserves the Poisson structure $\piG \oplus \piG$. 
Denote the quotient Poisson variety by
$ (G^{1, w^{-1}w_0},\,\piG)  \times_H (G^{1, w_0}, \;\piG)/H$.
We can now state and prove the main result of $\S$\ref{sec-log-pi}.

\bth{th-sec5-main}
For any $t \in H$ and $w \in W$, 
\begin{align*}
\rho: &\;(G^{1, w^{-1}w_0}, \; \piG)  \times_H (G^{1, w_0} , \;\piG)/H \lra 
(X_{t, w}, \, \Pi): \\  
&[g_1, \; {g_2}_\cdot H] \longmapsto [g_1g_2\dw_0, \; t(g_2\dw_0)_+]
\end{align*}
is a biregular Poisson isomorphism from 
$(G^{1, w^{-1}w_0}, \, \piG)  \times_H (G^{1, w_0} , \,\piG)/H$ to the
Zariski open subset $X_{t, w} \cap X_t^0$ of $X_{t, w}$. Moreover,
\[
X_{t, w} \cap X_t^0 =\{[n (m\dw_0^{-1})_+\dw_0, \, tm ]: \, 
n \in U \cap (B_- w^{-1}w_0 B_-), m  \in U \cap (B_- w_0 B_-)\}.
\]
\eth

\begin{proof}
Recall the Poisson map $\sigma$ in \prref{pr-key}. Replace $X_{t, w_0}$ in
$\sigma$ by its open subvariety $X_{t, w_0}^1$. It follows from \leref{le-Jt-poi} that
$\rho$ is Poisson. Identify
\[
(U \cap (B_-w^{-1}w_0B_-)) \times (U \cap (B_-w_0B_-)) \stackrel{\cong}{\lra}
(G^{1, w^{-1}w_0} \times G^{1, w_0})/(H \times H)
\]
by $(n, m) \mapsto (n, m)_\cdot (H \times H)$ for
$n \in U \cap (B_-w^{-1}w_0B_-)$ and $m \in U \cap (B_-w_0B_-)$, so that
$\rho$ becomes 
\[
{\rho}^\prime:  \; (U \cap (B_-w^{-1}w_0B_-)) \times (U \cap (B_-w_0B_-))\lra
X_{t, w}: \; (n, m) \longmapsto [nm\dw_0, \, t(m\dw_0)_+].
\]
Let 
$\xi: U \cap (B_- w_0 B_-) \to U \cap (B_- w_0 B_-): m \mapsto (m\dw_0)_+$ be
the isomorphism from \leref{le-xi}. 
The composition 
$\rho^{\prime\prime} := \rho^\prime ({\rm id} \times \xi^{-1})$ is then given by
\[
\rho^{\prime\prime}: \; (U \cap (B_-w^{-1}w_0B_-)) \times (U \cap (B_-w_0B_-))\lra
X_{t, w}: \; (n, m) \longmapsto [n\xi^{-1}(m)\dw_0, \, tm].
\]
By \leref{le-xi}, the image of $\rho^{\prime\prime}$ is given by
\[
{\rm Im} \rho^{\prime\prime} = 
\{[n (m\dw_0^{-1})_+\dw_0, \, tm ]: \, 
n \in U \cap (B_- w^{-1}w_0 B_-), m  \in U \cap (B_- w_0 B_-)\}.
\]
To prove the theorem, it now suffices to show that $\rho^{\prime\prime}$ is 
injective with image $X_{t, w} \cap X_t^0$, using the fact that $X_{t,w}$
is smooth.
Note that ${\rm Im} \rho^{\prime\prime}$ lies in the affine chart 
$X_t^\prime :=\{[n\dw_0, tm]: n, m \in U\}$
of $X_t$. 
It is clear from the explicit formula of $\rho^{\prime\prime}$ that it is
injective.
It remains to show that
$X_{t, w} \cap X_t^0 ={\rm Im} \rho^{\prime\prime}$. Clearly,
$X_{t, w} \cap X_t^0 \supset {\rm Im} \rho^{\prime\prime}$. 
Suppose that $[n_1 \dw_0, tn_2] \in X_{t, w} \cap X_t^0$. Then
$n_2 \in U \cap (B_- w_0 B_-)$ and $n_1 \dw_0 t n_2 \dw_0^{-1} n_1^{-1}
\in BwB_-$. Thus
\[
(n_2 \dw_0^{-1})_+ n_1^{-1} \in B_- w_0 B w B_- = B_- w_0 w B_-.
\]
Let $n = n_1 (n_2 \dw_0^{-1})_+^{-1}$. Then $n \in U \cap (B_- w^{-1}w_0 B_-)$
and $n_1 = n (n_2 \dw_0^{-1})_+$. Then $[n_1 \dw_0, tn_2] 
 = \rho^{\prime\prime}(n,n_2) \in {\rm Im} \rho^{\prime\prime}$.
\end{proof}

Recall that $\mu: G \times_B B \to G$ is the Grothendieck map.

\bco{co-sec5-main}
Let $G$ be simply connected. For any $t \in H$ and $w \in W$, the map
\begin{align*}
\mu {\rho}: &\;\;(G^{1, w^{-1}w_0}, \; \piG)  \times_H (G^{1, w_0} , \;\piG)/H \lra 
(F_{t, w}, \, \pi): \\
&\;\;[g_1, \; {g_2}_\cdot H]\longmapsto
g_1 g_2 \dw_0 t (g_2 \dw_0)_+ (g_2 \dw_0)^{-1} g_1^{-1}
\end{align*}
is a biregular Poisson isomorphism from 
$(G^{1, w^{-1}w_0}, \,\piG)  \times_H (G^{1, w_0} , \,\piG)/H$ to a
Zariski open subset of $(F_{t, w}, \pi)$.
\eco

\bre{re-logcanonical} 
In a future paper, we will use the rational isomorphisms $\rho$ and $\mu \rho $
to study log-canonical coordinates for the Poisson varieties $(X_{t,w}, \Pi)$
and $(F_{t, w}, \pi)$  and study
the associated cluster algebras.
\ere

\sectionnew{Appendix}\lb{appendix}

\subsection{Poisson Lie groups}\lb{subsec-poi-lie}
In this appendix, we recall some general facts on Poisson Lie groups
that are used in the construction of the Poisson structure 
$\pi$ on $G$ in $\S$\ref{sec-pi0-on-G} and the Poisson structure $\Pi$ on $G \times_B B$ in
$\S$\ref{sec-pi}.
Some of the omitted details in this section can be found in \cite{anton-malkin}
and \cite{k-s:quantum}.

Recall that a Poisson bi-vector field $\piG$ on a Lie group $G$ is said to be
{\it multiplicative} if the map 
$m: G \times G \to G: (g_1, g_2) \mapsto g_1g_2$ is Poisson with respect to
$\piG$. A {\it Poisson Lie group} is a Lie group $G$ with a multiplicative Poisson 
bi-vector field $\piG$. 
An action  $\sigma: G \times P \to P$
of a  Poisson Lie group $(G, \piG)$ 
on a Poisson manifold $(P, \pi_{{}_P})$ is said to be Poisson if 
$\sigma$ is  a Poisson map. 

If $(G, \piG)$ is a Poisson Lie group, then $\piG(e) = 0$, where $e \in G$ is the identity element.
The linearization of $\piG$ at $e$ is the linear map 
$\delta_\g: \g \to \wedge^2\g$ given by $\delta_\g(x) = [\tilde{x}, \piG](e)$, where
for $x \in \g$, $\tilde{x}$ is any vector field on $G$ with $\tilde{x}(e) = x$, and 
$[\tilde{x}, \piG]$ is the Lie derivative of $\piG$ by $\tilde{x}$.

Two Poisson Lie groups $(G, \piG)$ and $(G^*, \piGs)$ are said to be {\it dual} to
each other if their Lie algebras $\g$ and $\g^*$ are dual to each other and if
the dual map of $\delta_\g: \g \to \wedge^2 \g$ is the Lie bracket on $\g^*$ 
and the dual map of 
$\delta_{\g^*}: \g^* \to \wedge^2 \g^*$
is the Lie bracket on $\g$. 

One important class of Poisson Lie groups is constructed from {\it Manin triples}.
Recall that a Manin triple is a quadruple $(\d, \g, \g^*, \lara)$,
where $\d$ is an even dimensional Lie algebra, 
$\lara$ is a symmetric non-degenerate invariant
bilinear form on $\d$, $\g$ and $\g^*$ are Lie subalgebras of
$\d$, both maximally isotropic with respect to $\lara$, and $\g \cap \g^*
=0$. The bilinear form $\lara$ gives rise to a non-degenerate
pairing between $\g$ and $\g^*$, so $\g^*$ can indeed be regarded
as the dual space of $\g$.

Assume that $(\d, \g, \g^*, \lara)$ is a Manin triple. 
Let $\{x_i\}$ be a basis of $\g$ and let $\{\xi_i\}$ be
the dual basis of $\g^*$. Let
\[
\Lambda = \frac{1}{2} \sum_i (\xi_i \wedge x_i) \in \wedge^2 \d.
\]
Then $\Lambda$ is independent of the choices of the bases,
and the Schouten bracket $[\Lambda, \Lambda] \in \wedge^3\d$ is ad-invariant.
Let $D$ be a connected Lie group with Lie algebra $\d$.
Define
the bi-vector fields $\pi_{{\scriptscriptstyle D}}^{\pm}$ on $D$ by
\[
\pi_{{\scriptscriptstyle D}}^{\pm} = \Lambda^R \pm \Lambda^L,
\]
where $\Lambda^R$ and $\Lambda^L$ are respectively the right and left invariant
bi-vector fields on $D$ with $\Lambda^R(e) = \Lambda^L(e) = \Lambda$.
Then both $\pi_{{\scriptscriptstyle D}}^{-}$ and $\pi_{{\scriptscriptstyle D}}^+$
are Poisson structures on $D$ (see \cite[Proposition 3.4.1]{k-s:quantum}
for a proof that $\pi_{{\scriptscriptstyle D}}^{-}$ is Poisson, and use
the fact that left and right-invariant vector fields commute to
see $[\pi_{{\scriptscriptstyle D}}^+, \pi_{{\scriptscriptstyle D}}^+]
= [\pi_{{\scriptscriptstyle D}}^-, \pi_{{\scriptscriptstyle D}}^-]$,
so $\pi_{{\scriptscriptstyle D}}^+$ is Poisson).
Let $G$ and $G^*$
be the connected subgroups of $D$ with Lie algebras $\g$ and $\g^*$ 
respectively. 
One then checks that both $G$ and $G^*$ are Poisson submanifolds
of $(D, \pi_{{\scriptscriptstyle D}}^{-})$. Set
\[
\piG = \pi_{{\scriptscriptstyle D}}^{-}|_{{}_G} \hs \mbox{and} \hs 
\piGs = -\pi_{{\scriptscriptstyle D}}^{-}|_{{}_{G^{*}}}.
\]
Then 
$(G, \piG)$ and $(G^*, \piGs)$ is a pair of dual Poisson Lie groups,
and $(G, \piG)$ and $(G^*, -\piGs)$ are Poisson subgroups
 of $(D, \pi_{{\scriptscriptstyle D}}^{-})$.
In particular,
every Poisson action of $(D, \pi_{{\scriptscriptstyle D}}^{-})$ on 
a Poisson manifold 
restricts to Poisson actions of $(G, \piG)$ and $(G^*, -\piGs)$. 

The following \leref{le-D-pipm} is immediate from definitions.

\ble{le-D-pipm}
The following two actions are Poisson:
\begin{align*}
&(D, \pi_{{\scriptscriptstyle D}}^{-}) \times (D, \pi_{{\scriptscriptstyle D}}^{+}) 
\lra (D, \pi_{{\scriptscriptstyle D}}^{+}): \; \; (d_1, d_2) \longmapsto d_1d_2 \\
&(D, \pi_{{\scriptscriptstyle D}}^{+}) \times (D, -\pi_{{\scriptscriptstyle D}}^{-}) 
\lra (D, \pi_{{\scriptscriptstyle D}}^{+}): \; \; (d_1, d_2) \longmapsto d_1d_2.
\end{align*}
\ele

A proof of the following \prref{pr-DtoDG} can be found in \cite{anton-malkin}.
See also \cite{brown-goodearl-yakimov} and \cite{lu-yakimov-3}.
We give an outline of the proof for completeness.

\bpr{pr-DtoDG}
Assume that $G$ is closed in $D$, and let $p: D \to D/G$ be the natural projection. Then
\[
p(\pi_{{\scriptscriptstyle D}}^{\pm}) = p (\Lambda^R)
\]
is a well-defined Poisson structure on $D/G$, and the actions
\begin{align*}
& (G, \piG) \times (D/G, \; p(\pi_{{\scriptscriptstyle D}}^\pm)) \lra 
(D/G, \;p(\pi_{{\scriptscriptstyle D}}^\pm)): \; \; (g, d_\cdot G) \longmapsto 
gd_\cdot G,\\
& (G^*, -\piGs) \times (D/G, \; p(\pi_{{\scriptscriptstyle D}}^\pm)) 
\lra (D/G, \;p(\pi_{{\scriptscriptstyle D}}^\pm)): \; \; (u, d_\cdot G) \longmapsto 
ud_\cdot G
\end{align*}
are Poisson. Moreover, symplectic leaves of $p(\pi_{{\scriptscriptstyle D}}^\pm)$ 
in $D/G$ are precisely the connected components
of the nonempty intersections of $G$ and $G^*$-orbits in $D/G$.
\epr

\begin{proof}
The element $\Lambda \in \wedge^2 \d$ is mapped to $0$ under 
the projection $\d \to \d/\g$. Thus $p (\Lambda^L) = 0$. 
Since $p(\Lambda^R)$ is clearly well-defined,
$p(\pi_{{\scriptscriptstyle D}}^+) = p(\pi_{{\scriptscriptstyle D}}^-)$
 are well-defined. 
Now the Poisson action of
$(D, \pi_{{\scriptscriptstyle D}}^{-})$ on $(D, \pi_{{\scriptscriptstyle D}}^{+})$ 
by left multiplication 
restricts to give Poisson actions of $(G, \piG)$ and $(G^*, -\piGs)$, which
clearly descend to give Poisson actions on $(D/G, \;p(\pi_{{\scriptscriptstyle D}}^\pm))$.

Since $\g + \gst = \d$, 
the $G$-orbit $G_\cdot \ud$ and the $G^*$-orbit
$G^*{}_\cdot \ud$ intersect transversally at $\ud=d_\cdot G$ for any $d \in D$. In particular,
\[
T_{\ud} ((G_\cdot \ud) \cap (G^*{}_\cdot \ud)) 
= T_{\ud}(G_\cdot \ud) \cap T_{\ud}(G^*{}_\cdot \ud). 
\] 
To prove the statement about the symplectic leaves of $p(\pi_{{\scriptscriptstyle D}}^+)$, 
it is thus enough to check that
$T_{\ud}(G_\cdot \ud) \cap T_{\ud}(G^*{}_\cdot \ud)$ 
coincides with the tangent
space to the symplectic leaf of $p(\pi_{{\scriptscriptstyle D}}^+)$ at $\ud$. 
To this end, identify
$T_{\ud} (D/G) \cong \d/\Ad_d \g.$
Then $p(\pi_{{\scriptscriptstyle D}}^+)(\ud)$ 
becomes the element $p_d(\Lambda) \in \wedge^2 (\d/\Ad_d \g)$, where
$p_d: \d \to \d/\Ad_d \g$ is the projection. Let $\tilde{\Lambda}: \d \to \d$ be the
map
\[
\tilde{\Lambda}(x+\xi) = \frac{1}{2} \sum_{i=1}^{n} (\la x+\xi, \; \xi_i\ra x_i -
\la x + \xi, \; x_i \ra \xi_i) = \frac{1}{2}(x-\xi), \hs x \in \g, \xi \in \g^*.
\]
Let $S_{\ud}$ be the symplectic leaf of $p(\pi_{{\scriptscriptstyle D}}^+)$ through $\ud$. 
By definition,
\[
T_{\ud} S_{\ud} = p_d \left(\tilde{\Lambda}(\Ad_d \g)\right) \subset \d/\Ad_d \g.
\]
Since for any $x + \xi \in \Ad_d \g$, 
\begin{equation}\lb{eq-appendixlambda}
\tilde{\Lambda}(x+\xi) = \frac{1}{2}(x-\xi) =x - \frac{1}{2}(x+\xi) = -\xi + 
\frac{1}{2}(x+\xi),
\end{equation}
one sees that $p_d(\tilde{\Lambda}(\Ad_d \g)) = p_d (\g) \cap p_d (\g^*) \cong
T_{\ud}(G_\cdot \ud) \cap T_{\ud}(G^*{}_\cdot \ud)$. 
\end{proof}

\ble{le-Gs-pi0}
The local diffeomorphism $p|_{G^*}: (G^*, -\piGs) \to 
(D/G, p(\pi_{{\scriptscriptstyle D}}^+)): u \mapsto
u_\cdot G$ is Poisson.
\ele

\begin{proof} This is because $(G^*, -\piGs)$ is a Poisson subgroup of $(D, 
\pi_{\scriptscriptstyle D}^-)$.
\end{proof}

The following \leref{le-leaves-piD} on the symplectic leaves of 
$\pi_{{\scriptscriptstyle D}}^+$ in $D$ 
 is proved in \cite{anton-malkin}.
We give a slightly different proof here for completeness. Moreover, 
\eqref{eq-tilde-piD-1} in our
 proof of \leref{le-leaves-piD} 
is used in the proof of \prref{pr-Qt-sub}.

\ble{le-leaves-piD}
Symplectic leaves of $\pi_{{\scriptscriptstyle D}}^{+}$ in $D$ are precisely the connected components of
nonempty intersections of $(G, G^*)$-double cosets and $(G^*, G)$-double
cosets in $D$.
\ele

\begin{proof}
Let $d \in D$. Since
\[
T_d(G^*dG) + T_d(GdG^*) = r_d \g^* + l_d \g + r_d \g + l_d \g^* = r_d \d = T_d D,
\]
the two cosets $G^*dG$ and
$GdG^*$ intersect transversally at $d$. Let $\Sigma$ be the 
symplectic leaf of $\pi_{{\scriptscriptstyle D}}^{+}$ through $d$. It is enough to show that
\[
T_d\Sigma = T_d (G^*dG) \cap T_d (GdG^*).
\]
Let $\tilde{\pi}_{{\scriptscriptstyle D}}^{+}: T^*D \to TD$ be the bundle map
defined by $\pi_{{\scriptscriptstyle D}}^{+}$ (see \eqref{eq-tilde-pi}).
Identify $\d$ with $\d^*$ via $\lara$, and for $x \in \g$ and $\xi \in \g^*$, let $\alpha_{x+\xi}$
be the right invariant $1$-form on $D$ with value $x + \xi$ at the identity element of $D$.
By \eqref{eq-appendixlambda},
\begin{align}
\tilde{\pi}_{{\scriptscriptstyle D}}^{+} (\alpha_{x+\xi})(d) &= 
r_d \tilde{\Lambda} (x+\xi) + l_d \tilde{\Lambda} (\Ad_{d}^{-1}(x+\xi))\\
\lb{eq-tilde-piD-1} &=-r_{d} (\xi) + 
l_d \left(p_\g \Ad_{d}^{-1}(x+\xi)\right)\\
\lb{eq-tilde-piD-2}
&= r_d (x) - l_d \left(p_{\g^*} \Ad_{d}^{-1}(x+\xi)\right),
\end{align}
where $p_\g: \d \to \g$ and $p_{\g^*}: \d \to \g^*$ are the projections 
with respect to the decomposition $\d = \g + \g^*$.
Thus $T_d \Sigma \subset T_d (G^*dG) \cap T_d (GdG^*)$. Conversely, if
 $v_d \in T_d (G^*dG) \cap T_d (GdG^*)$, then
\[
v_d = -r_d (\xi) + l_d (x_1) = r_d(x) - l_d (\xi_1)
\]
for some $x, x_1 \in \g$ and $\xi, \xi_1 \in \g^*$, so $x_1 + \xi_1 = 
\Ad_{d}^{-1}(x+\xi)$, so
 $v_d = \tilde{\pi}_{{\scriptscriptstyle D}}^{+}(\alpha_{x+\xi})(d) 
\in T_d \Sigma$.
Hence $T_d \Sigma = T_d (G^*dG) \cap T_d (GdG^*)$.

\end{proof}
 
\subsection{Coisotropic reduction}\lb{subsec-coiso-reduction}

In this section, we prove \prref{pr-appendix-sub},
which is used in the study of the
Poisson structure $\pi$ on $G \times_B B$ in $\S$\ref{sec-pi}. \prref{pr-appendix-sub}
can be proved by combining Corollary 2.3 in \cite{marsden-ratiu-reduction}
and several examples therein, but we give a direct proof here for completeness.

A {\it Poisson vector space} is by definition a 
pair $(V, \pi)$, where $V$ is a vector space and $\pi \in 
\wedge^2 V$. Let $(V, \pi)$ be a finite dimensional 
Poisson vector space. Let
\[
\tilde{\pi}: \; \; V^* \lra V: \; \; (\tilde{\pi}(\xi), \; \eta) = \pi(\xi, \eta), 
\hs \xi, \eta \in V^*
\]
and set $V_\pi = \tilde{\pi}(V^*) \subset V$. A subspace
$V_1$ of $V$ is called a Poisson subspace of $(V, \pi)$
if $V_1 \supset V_\pi$, or, equivalently, if
$\pi \in \wedge^2 V_1$. In this case, $(V_1,\pi)$ is a Poisson vector
space.
Recall \cite{we:coiso} 
that a subspace $U$ of $V$ is said to be coisotropic with 
respect to $\pi$ if
$\tilde{\pi}(U^0) \subset U$, where
$U^0 = \{\xi \in V^*: \; \; (\xi, x) = 0, \; \forall x \in U\}.$

\bre{re-coisotropic} It is easy to see that 
$U \subset V$ is coisotropic if and only if 
$\pi$ is in the subspace $ U \wedge V$ of $\wedge^2 V$.
\ere



\ble{le-U-V}
Let $U$ be a coisotropic subspace of  $(V, \pi)$, let
$\phi: V \to V/\tilde{\pi}(U^0)$ be the projection, and set
$\varpi = \phi(\pi)$. Then
\begin{equation}\lb{eq-varpi}
(V/\tilde{\pi}(U^0))_\varpi = \phi(U \cap V_\pi) \subset U/\tilde{\pi}(U^0)\subset
 V/\tilde{\pi}(U^0).
\end{equation}
In particular, $(U/\tilde{\pi}(U^0), \, \varpi)$ is a Poisson vector subspace
of $(V/\tilde{\pi}(U^0), \, \varpi)$.
\ele

\begin{proof}
By a linear algebra computation, 
 $ (\tilde{\pi}(U^0))^0 = \tilde{\pi}^{-1}(U)=
\{\xi \in V^*: \tilde{\pi}(\xi) \in U\}$. By identifying
$(V/\tilde{\pi}(U^0))^* =  (\tilde{\pi}(U^0))^0$, one has 
\[
(V/\tilde{\pi}(U^0))_\varpi = \tilde{\varpi}((V/\tilde{\pi}(U^0))^*)
= \phi(\tilde{\pi}(\tilde{\pi}^{-1}(U))) = \phi(U \cap V_\pi).
\]
\end{proof}

Let $(P, \piP)$  be a Poisson manifold. A submanifold $Q \subset P$ is
 said to be
coisotropic if $\tpiP((T_qQ)^0) \subset T_qQ$ for every $q \in Q$, where
$(T_qQ)^0$ 
is the conormal bundle of $Q$ in $P$ 
and $\tpiP$ is the bundle map from $T^*P$ to $TP$ given in \eqref{eq-tilde-pi}. 
In this case, $\tpiP((T_qQ)^0)$ is called the {\it characteristic 
distribution} of $\piP$ on $Q$.


\bpr{pr-appendix-sub}
Let $(P, \piP)$ and $(R, \piR)$ be two Poisson manifolds with a 
surjective Poisson submersion $\phi: (P, \piP) \to
(R, \piR)$.
Assume that

1) $Q$ is a coisotropic submanifold 
of $(P, \piP)$ such that the characteristic distribution of $\piP$ on $Q$
coincides with the distribution defined by the tangent spaces to the fibers of $\phi$
and that $\phi(Q)$ is smooth submanifold of $R$; 

2) for every $q \in Q$, $Q$ intersects with the symplectic leaf $S_q$ of $\piP$ 
at $q$ cleanly, i.e.,
$Q \cap S_q$ is a submanifold and $T_q(Q \cap S_q) = T_q Q \cap T_qS_q$.

Then $\phi(Q)$ is a Poisson submanifold of
$(R, \piR)$, and for each $q \in Q$,
the symplectic leaf of $\piR$ in $\phi(Q)$ at
$\phi(q)$ is the connected component of $\phi(Q \cap S_q)$ containing $\phi(q)$. 
\epr

\begin{proof}
To show $\phi(Q)$ is a Poisson submanifold, it suffices to show
$T_{\phi(q)}(\phi(Q))$ is a Poisson subspace of $T_{\phi(q)}(R)$ for each
$q\in Q$. This last assertion is a consequence of the last statement
 of \leref{le-U-V} and assumption 1), since $Q$ is coisotropic.
 Furthermore, by 2), \eqref{eq-varpi} in 
\leref{le-U-V} gives the assertion on symplectic leaves.
\end{proof}

\subsection{Singularities of intersections of Bruhat cells and
Steinberg fibers}\lb{subsec-singsteinberg}

For an affine variety $X$ with ring of regular functions $O(X)$
and $g_1, \dots, g_k \in O(X)$, let $V(g_1, \dots, g_k)$ denote
the common vanishing set of $g_1, \dots, g_k$, and let
$(g_1,\dots, g_k)$ denote the ideal in $O(X)$ generated by
$g_1, \dots, g_k$. If $Y \subset X$ is Zariski closed, let
$I(Y)$ be the ideal of regular functions vanishing on $Y$.

Assume $G$ is simply connected for the remainder of the section.
Let $r = \dim H$, let $\{\alpha_1, \alpha_2, \ldots, \alpha_r\}$ be the set of simple roots, and 
let $\omega_1, \omega_2, \ldots, \omega_r$ be the corresponding {\it fundamental weights}, i.e., $\omega_j \in
\h^*$ for each $1\leq j \leq r$ and $\omega_j(h_{\alpha_k}) = \delta_{jk}$ for $1 \leq j, k \leq r$.
Denote by $\chi_j$ the character of the irreducible representation with $\omega_j$ as the highest weight.
Then the {\it Steinberg map} is the map \cite{hum-conjugacy}
\begin{equation}\lb{eq-chi}
\chi: \; G \lra \C^r: \; \; \chi(g) = (\chi_1(g), \, \chi_2(g), \ldots, \chi_r(g)).
\end{equation}
For $z = (z_1, \dots, z_r) \in \C^r$, let $F_z = F_t$ for any 
$t\in H$ such that $\chi(t)=z$.
For a Bruhat variety $\overline{BwB_-}$, let $f_i = {\chi_i}|_{\overline{BwB_-}} - z_i.$

\bpr{pr:smoothirredclosure}

\noindent 1) $F_z \cap \overline{BwB_-} = V(f_1, \dots, f_r)$.

\noindent 2) $F_z \cap \overline{BwB_-}$ is nonempty and irreducible

\noindent 3) $\dim(F_z \cap \overline{BwB_-})= d = \dim(G) - r - l(w)$.
\epr

\begin{proof}
Since $F_z = V(\chi_1 - z_1, \dots, \chi_r - z_r)$, part 1) is clear.
By \prref{pr-inter-Bruhat-Steinberg}, $F_z \cap \overline{BwB_-}$
is nonempty. Let $V$ be an
irreducible component of $V(f_1,\dots, f_r)$ and note that
 $\dim(V) \ge  \dim(G) - r - l(w)$.
Let $F_z = \cup_{i=1}^n C_{z_i}$
 be the decomposition of $F_z$ into
conjugacy classes with $C_{z_1}=R_z$ the unique regular conjugacy
class in $F_z$. Then
$ F_z \cap \overline{BwB_-} = \cup_{i=1, \dots, n, y \ge w } 
C_{z_i} \cap ByB_-$.
 If $i > 1$ or $y > w$, then
by \prref{pr-smoothirredint},  
\[
\dim(C_{z_i} \cap ByB_-) = \dim(C_{z_i}) - l(y) < \dim(R_z) - l(w) = d.
\]
It follows from \leref{le-irred} that $V(f_1, \dots, f_r)$ is irreducible, and
part 3) is an easy consequence.
\end{proof}

\ble{le-bruhatcm} 
$\overline{BwB_-}$ is Cohen-Macaulay for all $w\in W$.
\ele

\begin{proof}
By a Theorem of Ramanathan, 
$\overline{B_-wB_-}/B_-$ is  Cohen-Macaulay in $G/B_-$ 
 \cite{ra:cm}. 
The result now follows easily by using the fact that the smooth morphism
 $G \to G/B_-$ is a locally trivial bundle in
the Zariski topology and using the isomorphism
 $\overline{B_-wB_-} \cong \overline{Bw_0wB_-}$ given by
left multiplication by $\dw_0$.
\end{proof}

\ble{le-jantzen} (see \cite[Lemma 7.1]{ja:nilporbits})
Let $Y$ be an irreducible Cohen-Macaulay affine variety
of dimension $n$, and let $f_1, \dots, f_r \in O(Y)$.
Let $X=V(f_1, \dots, f_r)$. Suppose 

\noindent 1) $X$ is irreducible and

\noindent 2) there is a smooth point $y\in X$ such that
$df_1(y), \dots, df_r(y)$ are linearly independent.

Then $\dim (X) = n - r$ and the ideal 
$(f_1, \dots, f_r) = I(X)$.
\ele

\bre{re-jantzenproof}
Our statement is more general than the statement in \cite{ja:nilporbits}.
The proof is identical, once we recall a basic fact about
Cohen-Macaulay varieties. The ideal $(f_1, \dots, f_r) = Q_1 \cap \dots \cap
Q_s$ has a minimal primary decomposition.
The Cohen-Macaulay condition ensures that if 
$P_i = \sqrt{Q_i}$ for $i=1, \dots, s$,  the varieties
$V(P_i)$ all have the same dimension (by 
\cite[Theorem 17.6]{ma:ca}).
\ere

\bth{th-completeintersection}
The ideal $(f_1,\dots, f_r)$ is the ideal of functions vanishing on 
the irreducible variety 
 $F_z \cap \overline{BwB_-}$ in $\overline{BwB_-}$.
Moreover, $F_z \cap \overline{BwB_-}$
is Cohen-Macaulay.
\eth

\begin{proof}
By  \leref{le-bruhatcm}, $\overline{BwB_-}$ is Cohen-Macaulay.
By \prref{pr:smoothirredclosure}, $F_z \cap \overline{BwB_-}$
is irreducible. 
Recall that $(\chi_1 - z_1, \dots, \chi_r - z_r)$ is the ideal of
$F_z$ and $R_z$ is the smooth locus of $F_z$ (\cite{hum-conjugacy},
Theorem 4.24).
In particular,
\[
R_z = \{ y\in F_z : d\chi_1(y), \dots, d\chi_r(y) 
\ {\rm are \ linearly \ independent \ on  } \ T_y(G) \}
\]
and $T_y(R_z)$ is defined in $T_y(G)$ by the vanishing
of $d\chi_1(y), \dots, d\chi_r(y)$. 

Since $R_z$ and $BwB_-$ are smooth locally closed
subvarieties of $G$ and $R_z$ meets $BwB_-$ transversally, 
$R_z \cap BwB_-$ is smooth and locally closed in $F_z \cap \overline{BwB_-}$.
Let $y \in R_z \cap BwB_-$ and let $\lambda$ be a nonzero covector in
the span of $d\chi_1(y), \dots, d\chi_r(y)$. 
 Since $T_y(R_z) + T_y(BwB_-)
= T_y(G)$, it follows that
$\lambda$ is nonzero on $T_y(BwB_-)$. Thus, the restrictions
 $df_1(y), \dots, df_r(y)$ of $d\chi_1(y), \dots, d\chi_r(y)$
to $T_y(BwB_-)$ are
linearly independent. Since $T_y(BwB_-)=T_y( \overline{BwB_-})$,
we can
apply \leref{le-jantzen} to deduce the first assertion.
Since $\overline{BwB_-}$ is Cohen-Macaulay, 
 $F_z \cap \overline{BwB_-}$ is Cohen-Macaulay
using \cite[Corollary, page 65]{se:localg}.
\end{proof}

\bpr{pr-smoothlocus}

\noindent 1) Let ${\overline{BwB_-}}_{\rm ns} $ be the smooth locus of
$\overline{BwB_-}$. Then $R_z \cap {\overline{BwB_-}}_{\rm ns}$
is the smooth locus of $F_z \cap {\overline{BwB_-}}_{\rm ns}$.

\noindent 2) The singular locus of $F_z \cap \overline{BwB_-}$
has codimension at least $2$.

\epr

\begin{proof}
By \thref{th-completeintersection}, the ideal sheaf of
$F_z \cap {\overline{BwB_-}}_{\rm ns}$ is generated by $f_1, \dots, f_r$.
By the Jacobian criterion,  the smooth locus
of $F_z \cap {\overline{BwB_-}}_{\rm ns}$ is the set of points 
$y \in F_z \cap {\overline{BwB_-}}_{\rm ns}$ where
$df_1(y), \dots, df_r(y)$ are linearly independent on
$T_y({\overline{BwB_-}}_{\rm ns})$.
Let $y \in R_z \cap {\overline{BwB_-}}_{\rm ns}$ be in $BvB_-$.
 Using transversality as in
the proof of \thref{th-completeintersection}, it follows that
$df_1(y), \dots, df_r(y)$ are linearly independent on
$T_y(BvB_-)$, so they are linearly independent on
$T_y({\overline{BwB_-}}_{\rm ns}) \supset T_y(BvB_-)$.
Conversely, if $y\in F_z - R_z$, 
$d\chi_1(y), \dots, d\chi_r(y)$ are linearly dependent in
$T_y(G)$. As a consequence,
their restrictions $df_1(y), \dots, df_r(y)$
 are linearly dependent
 on $T_y({\overline{BwB_-}}_{\rm ns})$, which gives 1).
For 2), note that if 
$y$ is a singular point of $F_z \cap \overline{BwB_-}$,
then either $y$ is a singular point of $\overline{BwB_-}$
or $y$ is a singular point of $F_z \cap {\overline{BwB_-}}_{\rm ns}$.
Since the singular set of $\overline{BwB_-}$ has codimension at least
two \cite{bgg},
it suffices to show that the singular set of
$F_z \cap {\overline{BwB_-}}_{\rm ns}$ has codimension at
least two.  Let
$F_z = \cup_{i=1}^n C_{z_i}$
 be the decomposition into conjugacy classes with $C_{z_1}=R_z$.
By 1), the singular set of $F_z \cap {\overline{BwB_-}}_{\rm ns}$
is contained in 
\[
\cup_{v\ge z, i\ge 2} C_{z_i} \cap BvB_-.
\]
By \cite[4.24]{hum-conjugacy}, if $i\ge 2$,
$\dim(C_{z_i}) \le \dim(R_z) - 2$. Since
 $\dim(C_{z_i} \cap BvB_-) = \dim(C_{z_i}) - l(v)$ by
\prref{pr-smoothirredint}, it follows that if $i \ge 2$,
\[
\dim(C_{z_i} \cap BvB_-) \le \dim(R_z \cap BvB_-) - 2 \le
\dim(F_z \cap \overline{BwB_-}) - 2.
\]
Part 2) follows.
\end{proof}

\bth{th-normality}
$F_z \cap \overline{BwB_-}$ is normal.
\eth

\begin{proof} 
Since $F_z \cap \overline{BwB_-}$
is Cohen-Macaulay by \thref{th-completeintersection},
 condition $S_2$ of Serre is satisfied (see \cite[page 183]{ma:ca}). 
Part 2) of \prref{pr-smoothlocus} is equivalent to condition
$R_1$ of Serre, so the theorem follows using
 Serre's normality criterion
 (\cite[Theorem 23.8]{ma:ca}).
\end{proof}


\begin{thebibliography}{99}

\bibitem{anton-malkin}
Alekseev, A. and Malkin, A., {\em Symplectic structures associated to Lie-Poisson groups}, 
Comm. Math. Phys. {\bf 162}(1) (1994), 147 - 173. 

\bibitem{beauville}
Beauville, A., {\em Symplectic singularities},  Invent. Math.
{\bf 139}(3) (2000), 541 - 549. 

\bibitem{bgg}
Bernstein, I.N., Gelfand, I.M., and Gelfand, S.I., {\em Schubert cells
and the cohomology of spaces G/P}, Russ. Math. Surveys {\bf 28}(3)
(1973), 3 - 26.

\bibitem{beren-z}
Berenstein, A. and Zelevinsky, A., {\em Tensor product multiplicities, canonical bases and totally
positive varieties}, Invent. Math. {\bf 143} (2001), 77 - 128.

\bibitem{brown-goodearl-yakimov}
Brown, K., Goodearl, K., and Yakimov, M., {\em  Poisson structures 
on affine spaces and flag varieties. I. Matrix affine Poisson space },
math.RT/0501109, to appear in Adv. Math.

\bibitem{chriss-ginzburg}
Chriss, N., and Ginzburg, V., {\em Representation theory and complex
geometry}, Birkhauser, 1997.


\bibitem{DKP}
De Concini, C., Kac, V. G., Procesi, C., {\em Some remarkable degenerations of 
quantum groups}.  Comm. Math. Phys.  {\bf 157}(2) (1993), 405 - 427.


\bibitem{deo}
Deodhar, V., {\em On some geometric aspects of Bruhat orderings, I, a finer
decomposition of Bruhat cells}, Invent. Math. {\bf 79}(3) (1985), 499 - 511.

\bibitem{drinfi}
Drinfeld, V., {\em Quantum groups},  Proceedings of the International Congress of 
Mathematicians, Vol. 1, 2 (Berkeley, Calif., 1986),  798 - 820, Amer. Math. Soc., 
Providence, RI, 1987.

\bibitem{ellers-gordeev} 
Ellers, E. and Gordeev, N., {\em Intersection of conjugacy classes with
Bruhat cells in Chevalley groups}, Pac. J. of Math. {\bf 214} (2) (2004),
245 - 261.

\bibitem{e-l:real}
Evens, S. and J.-H. Lu, {\em On the variety of Lagrangian subalgebras, I},
Ann. Ecole Norm. Sup. {\bf 34} (4)   (2001), 631 - 668.

\bibitem{e-l:cplx}
Evens, S. and J.-H. Lu, {\em On the variety of Lagrangian subalgebras, II},
Ann. Ecole Norm. Sup. {\bf 39} (2) (2006), 347 - 379.



\bibitem{FZ}
Fomin, S. and Zelevinsky, A., {\em Double Bruhat cells and total positivity}, 
J. Amer. Math. Soc., {\bf 12} (2) (1999), 335 - 380.

\bibitem{fz:I}
Fomin, S. and Zelevinsky, A.,
Cluster algebras. I. Foundations, {\em J. Amer. Math. Soc.}
{\bf 15} (2002), no. 2, 497 - 529 (electronic). 

\bibitem{fz:II}
Fomin, S. and Zelevinsky, A., 
Cluster algebras II: Finite type classification,
{\em  Invent. Math. } {\bf 154} (2003), no. 1, 63 - 121.

\bibitem{fu-symplectic}
Fu, B., {\em Symplectic resolutions for nilpotent orbits},  Invent. Math.,
{\bf 151} (1) (2003), 167 - 186.

\bibitem{fu-poisson} 
Fu, B., {\em Poisson resolutions},  J. Reine Angew. Math., {\bf 587}
 (2005), 17 - 26. 

\bibitem{fu-survey}
Fu, B., {\em A survey on symplectic singularities and resolutions},
Annales Mathématiques Blaise Pascal, {\bf 13} (2) (2006), 209 - 236.


\bibitem{goodearl-yakimov}
Goodearl, K. and Yakimov, M., {\em Poisson structures on affine spaces
and flag varieties, II}, math.QA/0509075.

\bibitem{HL}
Hodges, T. and Levasseur, T., {\em Primitive ideals of ${\mathbb C}_q[SL(3)]$}, 
Comm. Math. Phys. {\bf 156} (3) (1993), 581 - 605.

\bibitem{HKKR}
Hoffmann, T., Kellendonk, J., Kutz, N., and Reshetikhin, N., 
{\em Factorization dynamics and Coxter-Toda lattices}, Comm. Math. Phys. 
{\bf 212} (2) (2000), 297 - 321.

\bibitem{hum-conjugacy}
Humphreys, J., {\em Conjugacy classes in semisimple algebraic groups,} Mathematical 
surveys and monographs, Vol. 43, AMS, 1995.

\bibitem{ja:nilporbits}
Jantzen, J., {\em Nilpotent orbits in representation theory} in
{\em Lie Theory, Lie algebras and representations},  Birkhauser, 2004.

\bibitem{kasym}
Kaledin, D., {\em Symplectic resolutions: deformations and birational maps},
math.AG/0012008.

\bibitem{KL} Kazhdan, D. and Lusztig, G., {\em Representations of Coxeter groups 
and Hecke algebras}.  Invent. Math.  {\bf 53}(2) (1979), 165 - 184.

\bibitem{KZ}
Kogan, M. and Zelevinsky, A., {\em On symplectic leaves and 
integrable systems in standard complex semisimple Poisson-Lie
groups}, Int. Math. Res. Not. {\bf 32} (2002), 1685 - 1702. 

\bibitem{k-s:quantum}
Korogodski, L. and Soibelman, Y., {\em Algebras of functions on
quantum groups, part I}, AMS, Mathematical surveys and
monographs, Vol. 56, 1998.

\bibitem{camille}
Laurent-Gengoux, C., {\em From Lie groupoids to resolutions of singularity. Applications to symplectic resolutions. (I)}, math.DG/0610288.

\bibitem{lu-we:poi}
Lu, J.-H. and Weinstein, A., {\em Poisson Lie groups, dressing transformations,
and Bruhat decompositions},  J. Diff. Geom. {\bf 31} (2) (1990), 
501 - 526.  

\bibitem{lu-yakimov-3}
Lu, J.-H. and Yakimov, M., {\em 
Group orbits and regular partitions of Poisson manifolds}, math.SG/0609732.


\bibitem{lusztig}
Lusztig, G., {\em Total positivity in reductive groups}, {\em Lie theory and
geometry: in honor of Bertram Kostant}, 531 - 568, Progress in Math. {\bf 123}, Birkhauser, 1994.

\bibitem{marsden-ratiu-reduction}
Marsden, J. and Ratiu, T., {\em Reduction on Poisson manifolds}, Lett. Math. Phys.,{\bf 11} (1986), 161 - 169.

\bibitem{ma:ca}
Matsumura, H., {\em Commutative ring theory} Cambridge University Press,
1989.

\bibitem{ra:cm}
Ramanathan, A., {\em Schubert varieties are arithmetically Cohen-Macaulay},
Invent. Math. {\bf 80}, 283 - 294 (1985).


\bibitem{Reshe}
Reshetikhin, N., {\em Integrability of characteristic Hamiltonian systems on simple
Lie groups with standard Poisson Lie structure},  Comm. Math. Phys., {\bf 242}
(2003), 1 - 29.

\bibitem{Ri} R. Richardson, {\em{Intersections of double 
cosets in algebraic groups,}} Indagationes Mathematicae, Volume 3, Issue 1, (1992),
69 - 77.

\bibitem{se:localg}
Serre, J-P., {\em Local algebra}, Springer, 2000.



\bibitem{slodowy}
Slodowy, P., {\em Simple singularities and simple algebraic groups}, LNM
{\bf 815}, Springer, 1980.


\bibitem{steinberg-IHES}
Steinberg, R., {\em Regular elements of semisimple algebraic groups}, Inst. Hautes Etudes Sci. Publ.
Math. {\bf 25} (1965), 49 - 80.

\bibitem{steinberg-verlag}
Steinberg, R., {\em Conjugacy classes in algebraic groups} (notes by V.V. Deodhar), Lecture
Notes in Math. {\bf 366}, Springer-Verlag, Berlin, 1974.

\bibitem{steinberg-inventiones}
Steinberg, R., {\em On the desingularization of the unipotent variety}, Invent. Math. {\bf 36} (1976), 209 - 224.

\bibitem{we:coiso}
Weinstein, A., {\em Coisotropic calculus and Poisson groupoids}, 
 J. Math. Soc. Japan,
{\bf 40} (4) (1988), 705 - 726.

\end{thebibliography}
\end{document}